\documentclass[12pt,reqno]{amsart}
\usepackage{amsmath}
\usepackage{amsxtra}
\usepackage{amscd}
\usepackage{amsthm}
\usepackage{amsfonts}
\usepackage{amssymb}
\usepackage{eucal}

\newtheorem{theorem}{Theorem}[section]
\newtheorem{conj}[theorem]{Conjecture}
\newtheorem{cor}[theorem]{Corollary}
\newtheorem{lem}[theorem]{Lemma}
\newtheorem{prop}[theorem]{Proposition}

\theoremstyle{definition}

\theoremstyle{remark}

\theoremstyle{remark}

\numberwithin{equation}{section}

\newcommand{\nc}{\newcommand}
\nc{\on}{\operatorname}
\nc{\ch}{\mbox{ch}}
\nc{\Z}{{\mathbb Z}}
\nc{\C}{{\mathbb C}}
\nc{\pone}{{\mathbb C}{\mathbb P}^1}
\nc{\pa}{\partial}
\nc{\F}{{\mathcal F}}
\nc{\arr}{\rightarrow}
\nc{\larr}{\longrightarrow}
\nc{\al}{\alpha}
\nc{\ri}{\rangle}
\nc{\lef}{\langle}
\nc{\W}{{\mathcal W}}
\nc{\la}{\lambda}
\nc{\ep}{\epsilon}

\nc{\su}{\widehat{{\mathfrak sl}}_2}
\nc{\sw}{{\mathfrak s}{\mathfrak l}}

\nc{\g}{{\mathfrak g}}
\nc{\h}{{\mathfrak h}}
\nc{\n}{{\mathfrak n}}
\nc{\N}{\widehat{\n}}
\nc{\G}{\widehat{\g}}
\nc{\De}{\Delta_+}
\nc{\gt}{\widetilde{\g}}
\nc{\Ga}{\Gamma}
\nc{\one}{{\mathbf 1}}
\nc{\z}{{\mathfrak Z}}
\nc{\zz}{{\mathcal Z}}
\nc{\Hh}{{\mathcal H}_\beta}
\nc{\qp}{q^{\frac{k}{2}}}
\nc{\qm}{q^{-\frac{k}{2}}}
\nc{\La}{\Lambda}
\nc{\wt}{\widetilde}
\nc{\qn}{\frac{[m]_q^2}{[2m]_q}}
\nc{\cri}{_{\on{cr}}}
\nc{\kk}{h^\vee}
\nc{\sun}{\widehat{\sw}_N}
\nc{\hh}{\widehat{\mathfrak h}}
\nc{\HH}{{\mathcal H}_{q,t}}
\nc{\ca}{\wt{{\mathcal A}}_{h,k}(\sw_2)}
\nc{\gl}{\widehat{{\mathfrak g}{\mathfrak l}}_2}
\nc{\el}{\ell}
\nc{\s}{{\mathbf s}}
\nc{\bi}{\bibitem}
\nc{\om}{\omega}
\nc{\WW}{\W_\beta}
\nc{\scr}{{\mathbf S}}
\nc{\ab}{{\mathbf a}}
\nc{\rr}{r}
\nc{\ol}{\overline}
\nc{\con}{qt^{-1} + q^{-1}t}
\nc{\den}{q^{\el-1} t^{-\el+1}+ q^{-\el+1} t^{\el-1}}
\nc{\ds}{\displaystyle}
\nc{\B}{B}
\nc{\A}{{\mathbb A}}
\nc{\GG}{{\mathcal G}}
\nc{\UU}{{\mathcal U}}
\nc{\MM}{{\mathcal M}}
\nc{\CC}{{\mathcal C}}
\nc{\GL}{{}^L G}
\nc{\dzz}{\frac{dz}{z}}
\nc{\Res}{\on{Res}}
\nc{\rep}{{\mathcal R}ep \;}
\nc{\uqg}{U_q \G}
\nc{\uqgg}{U_q \g}
\nc{\Fq}{{\mathbb F}_q}

\nc{\stimes}{\ltimes}
\nc{\K}{\hat{\mathcal K}}
\nc{\Ql}{\ol{\mathbb Q}_\ell}

\nc{\ga}{\gamma}
\nc{\PL}{{}^L P}
\nc{\E}{\mc E}
\nc{\mc}{\mathcal}
\nc{\mbf}{\mathbf}
\nc{\bb}{{\mathfrak b}}
\nc{\OO}{{\mc O}}
\nc{\Po}{{\mc P}}
\nc{\V}{{\mc V}}
\nc{\yy}{{\mc Y}}
\nc{\M}{\mathcal M}
\nc{\Coh}{{{\mathcal C}oh}}
\nc{\Cohn}{\Coh_n}
\nc{\f}{{\mathcal F}}
\nc{\si}{_E}
\nc{\Gaf}{{\mathbb G}_{a,\Fq}}
\nc{\KK}{{\mathfrak k}}

\nc{\PCr}{{ \bs P  (\C[x])^r   }}
\nc{\PCN}{{ \bs P  (\C[x])^N   }}

\nc{\sN}{sl_{2N+1}}

\nc{\Pzr}{{ \bs P(\C((x-z)))^r}}

\nc{\PzN}{{ \bs P(\C((x-z)))^N}}

\newcommand{\bean}{\begin{eqnarray}}
\newcommand{\eean}{\end{eqnarray}}
\newcommand{\be}{\begin{displaymath}}
\newcommand{\ee}{\end{displaymath}}
\newcommand{\bea}{\begin{eqnarray*}}   
\newcommand{\eea}{\end{eqnarray*}}
\newcommand{\bs}{\boldsymbol}
\newcommand{\Ref}[1]{{$($\ref{#1}$)$}}

\begin{document}

\title[Critical points and flag varieties]{Critical points of master
  functions and flag varieties}
\author[E. Mukhin and A. Varchenko]
{E. Mukhin and A. Varchenko}
\thanks{Research of E.M. is supported in part by NSF grant DMS-0140460.
Research of A.V. is supported in part by NSF grant DMS-9801582}
\address{E.M.: Department of Mathematical Sciences, Indiana University -
Purdue University Indianapolis, 402 North Blackford St, Indianapolis,
IN 46202-3216, USA, \newline mukhin@math.iupui.edu}
\address{A.V.: Department of Mathematics, University of North Carolina 
at Chapel Hill, Chapel Hill, NC 27599-3250, USA, anv@email.unc.edu}
\maketitle

\begin{center} August 1, 2002
\end{center}

\begin{abstract}
We consider critical points of master functions associated 
with integral dominant weights
of Kac-Moody algebras and introduce a generating procedure constructing new critical points
starting from a given one. The set of all critical points constructed from a given one is
called a population. We 
formulate a conjecture that a population is isomorphic to the flag variety of the Langlands 
dual Kac-Moody algebra and prove the conjectures for algebras $sl_{N+1}, so_{2N+1}$, and $sp_{2N}$.

We show that populations associated with a collection of integral dominant $sl_{N+1}$-weights 
are in one to one correspondence with intersection points of suitable Schubert cycles in a
Grassmannian variety.
\end{abstract}

\section{Introduction}

Let $\h$ a the Cartan subalgebra of $sl_{r+1}$;
$( \ , \  )$  the Killing form on  $\h^*$;
$\al_1, \dots , \al_r \in \h^*$  simple roots;
$\La_1, \dots , \La_n \in \h^*$ dominant integral weights;
$l_1, \dots , l_r$ non-negative integers;
$z_1, \dots , z_n$ distinct complex numbers. The master function
$\Phi$ associated with this data 
is the rational function of $l_1 + \dots + l_r$ variables
$ t_1^{(1)}, \dots , t_{l_1}^{(1)}, \dots , t_{l_r}^{(r)}$,
given by  formula \Ref{master}.
We may think that  $l_1 + \dots + l_r + n$ particles are given in $\C$. The particle $t^{(i)}_j$
has weight $- \al_i$ and the particle $z_s$ has weight $\La_s$. The particles interact pairwise.
The interaction of particles $x$ and $y$ with weights $v$ and $w$, respectively, is given by
$(x - y)^{(v,w)}$. Then total interaction is the product of such terms
over all pairs. The master function describes the interaction of $t$-particles with
$t$-particles and $z$-particles. We ignore the interaction of $z$-particles with $z$-particles
since the positions of $z$ particles will be fixed in our considerations.

The master function appears in hypergeometric solutions to the KZ equations with values in
the tensor product of irreducible highest weight representations $L_{\La_1}, \dots , L_{\La_n}$
with highest weights $\La_1, \dots , \La_n$, respectively.
The solutions have the form
\be
u (z)\ = \ \int \ \Phi (\bs t, \bs z )^{{1\over \kappa}}\  A( \bs t; \bs z ) \ d \bs t\ ,
\ee
 where
$\kappa$ is the parameter of the KZ equations and
$A( \bs t; \bs z )$ is some explicitly written rational function with values in the tensor product
\cite{SV}. 

The master function also appears in the Bethe ansatz of the Gaudin model with values 
in the same tensor product \cite{RV}. 
In that case the value of the function $A(\ \cdot\ ; \bs z )$ at a point
$\bs t$ is an eigenvector of the commuting Gaudin Hamiltonians if $\bs t$ is a critical point
of the master function.

In this article we study critical points of the master function on the set where all $\{ t^{(i)}_j, z_s\}$
are distinct. In other words we study  those positions of distinct particles
$\{ t^{(i)}_j \}$ in the complement to $\{z_s\}$ which extremize the master function. 

The master function is symmetric with respect to permutations of variables  $\{ t^{(i)}_j \}$
with the same upper index. It is convenient to use symmetric coordinates.
To every point $\{ t^{(i)}_j \}$ we assign a tuple of polynomials $\bs y = (y_1, \dots , y_r)$ of one variable $x$
with $y_i(x) = \prod_j(x - t^{(i)}_j)$. The master function can be considered as a function on
the space  $\C_{l_1}^0[x] \times \dots \times \C_{l_r}^0[x]$, where $\C_d^0[x]$ is the 
affine space of monic polynomials of degree $d$.
We consider the critical set of the master function on this product.

The  Weyl group 
$\mathcal W$ acts on $\h^*$. For $w \in \mathcal W$ and $\la \in \h^*$ denote
$w\cdot\la=w(\la+\rho)-\rho$ the action of the Weyl group shifted by $\rho$, the half sum of positive roots.
Set $\La_\infty = \sum_{s=1}^n \La_s - \sum_{i=1}^r l_i \al_i$.

Let $B \subset SL(r+1, \C)$ be a Borel subgroup. The flag variety $SL(r+1, \C)/B$ has the Bruhat cell decomposition
$SL(r+1, \C)/ B\ =\ \sqcup_{w\in \mc W}\, G_w$ with $G_w \ =\  B w B / B$.

\begin{theorem}\label{th intr}
${}$

\begin{enumerate}
\item[$\bullet$] If $\La_\infty$ belongs to  one of the reflection hyperplanes of the shifted action
of the Weyl group, then the master function does not have critical points such that 
$\{ t^{(i)}_j, z_s\}$ are distinct.
\item[$\bullet$] 
Let $w_0 \in \mc W$ be such that $\tilde \La_\infty = w_0 \cdot \La_\infty$ is  dominant.
If $\tilde \La_\infty - \sum_{s=1}^n \La_s $ is not a linear combination of simple roots with
non-negative coefficients, then the master function does not have critical points such that 
$\{ t^{(i)}_j, z_s\}$ are distinct.
\item[$\bullet$]
Let $\tilde \La_\infty = w_0 \cdot \La_\infty$ be  dominant.
Assume that $\tilde \La_\infty - \sum_{s=1}^n \La_s $ is a linear combination of simple roots with
non-negative coefficients. 
Let a tuple $\bs y^0 \in \C_{l_1}^0[x] \times \dots \times \C_{l_r}^0[x]$ represent
a critical point of the master function such that  $\{ t^{(i)}_j, z_s\}$ are distinct. 
Then the closure in $\C_{l_1}^0[x] \times \dots \times \C_{l_r}^0[x]$
of the connected component of the critical set containing $\bs y^0$ is isomorphic to the Bruhat 
cell $G_{w_0}$.
\item[$\bullet$]
The number 
of connected components of the critical set is not greater, than the multiplicity of the representation
$L_{\tilde \La_\infty}$ in the tensor product $L_{\La_1} \otimes \dots \otimes L_{\La_n}$.
\end{enumerate}
\end{theorem}

Conjecturally 
for generic $z_1, \dots , z_n$ and arbitrary integral dominant weights
$L_{\La_1} \otimes \dots \otimes L_{\La_n}$
the number of connected components is equal to the multiplicity indicated in part 4 of the theorem.
For $sl_2$ the theorem and the conjecture are proved in \cite{ScV}. 

For $sl_2$ the 
theorem takes the following form. Irreducible representations are labeled by non-negative integers.
Denote $L_a$ the irreducible representation with highest weight $a \in \Z_{\geq 0}$, \ $\dim L_a = a + 1$.
The master function takes the form
\be
\Phi (\bs t; \bs z) \ = \ \prod_{j=1}^{l}\prod_{s=1}^n   (t_j - z_s)^{-m_s}
  \prod_{1\leq i < j \leq l} (t_i - t_j)^2 \ 
\ee
for given positive integers $m_1, \dots , m_n$. Set $|m| = m_1 + \dots + m_n$. To every $\bs t$ assign a polynomial
$y(x) = \prod_{j=1}^l (x - t_j)$ and consider the master function as a function on $\C_l^0[x]$.

\begin{theorem} \cite{ScV}
${}$

\begin{enumerate}
\item[$\bullet$] If $|m| + 1 - l < 0$ or $|m| + 1 - l = l$, then the master function does not have
critical points on the set where $\{ t_j, z_s\}$ are distinct.
\item[$\bullet$]
If $l < |m| + 1 - l $, then for generic $z_1, \dots , z_n$ the number of polynomials $\bs y$, such that
$\bs y$ represents a critical point and $\{ t_j, z_s \}$ are distinct, is equal to the multiplicity of 
the representation $L_{|m| - 2l}$ in the tensor product $L_{m_1}\otimes \dots \otimes L_{m_n}$. The critical
points are non-degenerate.
\item[$\bullet$]
Assume that $0 < |m| + 1 - l < l $, and  $z_1, \dots , z_n$ are generic.
 Let a tuple $\bs y^0 \in \C_{l}^0[x]$ represent a critical point of the master 
function such that  $\{ t_j, z_s\}$ are distinct. 
Then the closure in $\C_{l}^0[x]$ of the connected component of the critical set containing $\bs y^0$ is 
a straight line. The total number of the straight lines is equal to the
multiplicity of the representation
$L_{2l - 2- |m|}$ in the tensor product $L_{m_1}\otimes \dots \otimes L_{m_n}$.

\end{enumerate}
\end{theorem}

In Section 2  we consider master functions associated with a Kac-Moody algebra. In Section 3 
we describe a reproduction procedure which  constructs new critical points 
if a critical point is given. This is the first main ingredient of this paper.

We introduce the notion of the population of
critical points originated at a given critical point as the set of all 
critical points which are constructed from a given one by the reproduction procedure. 
Every population consists of
critical points of several master functions. We develop basic properties of populations and formulate
a conjecture that every population is isomorphic to the flag variety of the Kac-Moody 
algebra  Langlands dual
to the initial Kac-Moody algebra.
We  conjecture 
that the master functions whose critical points enter a given population are labeled by elements
of the Weyl group, and the part of the population, consisting of critical points of the 
master function labeled by an element $w \in \mc W$, is isomorphic to the Bruhat cell $G_w$ of the flag variety.
We conjecture that the number of populations originated at critical points of a given master function for generic
$z_1, \dots z_n$ is  given in terms of the representation theory
as  the multiplicity of a suitable irreducible representation in a suitable tensor product 
of irreducible representations, cf. Theorem \ref{th intr}.

In Section 4 we remind some statements from Schubert calculus. 

In Section 5 we collect
the results on $sl_{r+1}$-populations,  in particular  Theorem \ref{th intr}
is proved. 
In Section 5 we
 assign to every critical point $\bs y$ an $r+1$-dimensional space $V$
of polynomials in one variable.  This is the second main object of the paper.

We describe the 
ramification properties of $V$ (in the sense of Schubert calculus)  in terms
of complex numbers $z_1, \dots , z_n$,  
weights $\La_1, \dots , \La_n$, and the integers $l_1, \dots , l_n$, 
associated to the initial master function.
We prove that the population originated at $\bs y$ is naturally isomorphic to the full flag variety 
$FL(V)$ of the space $V$.
Knowing ramification 
properties of the spaces $V$, assigned to critical points of the initial master function,
and applying Schubert calculus we  give an estimate from above on the number of populations originated 
at critical points of the initial master function. 

We show that every $r+1$-dimensional vector space 
of polynomials of one variable corresponds to an  $sl_{r+1}$-population of critical points.

In Section 6 we develop  
the notion a selfdual vector space of polynomials of one variable. We show that a selfdual 
vector space has a canonical bilinear form $(\ , \ )$
which is symmetric if the dimension of the space is odd and is skew-symmetric
if the dimension is even. We discuss properties of the variety $FL^\perp(V) \subset FL(V)$ of full 
flags isotropic with
respect to the canonical bilinear form.

In  Section 7 we study $so_{2N+1}$ and $sp_{2N}$ populations.
We assign to every critical point $\bs y$ of a master function associated with $so_{2N+1}$ a selfdual
$2N$-dimensional space of polynomials $V$ in one variable. We describe the ramification properties of $V$ 
(in the sense of Schubert calculus)  in terms of numbers $\{z_s\}$,  weights $\{\La_s\}$ and the vector $\{l_s\}$ 
associated to the initial master function. 
We prove that the population originated at $\bs y$ is naturally isomorphic to the variety 
$FL^\perp(V)$ of isotropic flags of the space $V$.

We assign to every critical point $\bs y$ of a master function associated with $sp_{2N}$ a selfdual
$2N+1$-dimensional space of polynomials $V$ in one variable. We describe the ramification properties of $V$ 
(in the sense of Schubert calculus)  in terms of numbers $\{z_s\}$,  weights $\{\La_s\}$ and the vector $\{l_s\}$ 
associated to the initial master function. 
We prove that the population originated at $\bs y$ is naturally isomorphic to the variety 
$FL^\perp(V)$ of isotropic flags of the space $V$.

The notion of a selfdual space of polynomials of one variable and the
construction of the selfdual space of polynomials, associated to a critical point of
the $so_{2N+1}$ or $sp_{2N}$ type, is the third main object of the paper.

According to our conjecture there has to exist 
a Schubert calculus of selfdual $2N$-dimensional spaces of polynomials 
where the intersection numbers of
Schubert cycles would be given in terms of representation theory of $so_{2N+1}$. Analogously
there has to exist 
a Schubert calculus of selfdual $2N+1$-dimensional spaces of polynomials where 
the intersection numbers of
Schubert cycles would be given in terms of representation theory of $sp_{2N}$.

In Section 7 we formulate Theorems \ref{conv so} and \ref{conv sp} 
that every $2N$-dimensional selfdual space corresponds to an $so_{2N+1}$-population
and every $2N+1$-dimensional selfdual space corresponds to an $sp_{2N}$-population.

In Section 8 we sketch the theory of populations of germs of holomorphic functions of one variable
and then
prove Theorems \ref{conv so} and \ref{conv sp}. 

\newpage

In the Appendix useful facts about Wronskians are collected.

\medskip

The authors especially thank
 A. Oblomkov. In discussions with him we realized 
connections of the Schubert calculus with the problem of counting vector spaces of polynomials
of one variable with prescribed  ramification properties.
A. Oblomkov helped us to write the review Section 4.

We thank E. Frenkel. After we informed him about the results of the 
present work,E. Frenkel indicated to us connections of our work with
his research on Bethe ansatz and Miura
transformations \cite{FFR, Fr}. We plan to develop
those connections in a future joint work with E. Frenkel.

We thank A. Gabrielov, A. Eremenko and D. Novikov for useful discussions.

\section{Master functions and critical points}
\subsection{Kac-Moody algebras}\label{Kac_Moody sec}

Let $A=(a_{ij})_{i,j=1}^r$ be a generalized  Cartan matrix, 
$a_{ii}=2$,
$a_{ij}=0$ if and only $a_{ji}=0$,
 $a_{ij}\in \Z_{\leq 0}$ if $i\ne j$. 
We  assume that $A$ is symmetrizable, 
there is a diagonal matrix $D=\on{diag}\{d_1,\dots,d_r\}$ 
with positive integers $d_i$ such that $B=DA$
is symmetric. 

Let $\g=\g(A)$ be the corresponding complex Kac-Moody 
Lie algebra (see \cite{K}, \S 1.2), 
$\h \subset \g$  the Cartan subalgebra.
The associated scalar product is non-degenerate on $\h^*$ and 
 $\on{dim}\h = r + 2d$,  where $d$ is the dimension of 
the kernel of the Cartan matrix $A$.

Let $\al_i\in \h^*$, $\al_i^\vee\in \h$, $i = 1, \dots , r$, be the sets of simple roots,
coroots, respectively. We have
\bea
 (\al_i,\al_j)&=& d_i \ a_{ij}, \\
\langle\la ,\al^\vee_i\rangle&=&2(\la,\al_i)/{(\al_i,\al_i)}, \qquad \la\in\h^*.
\eea

Let $\mathcal P = \{ \lambda \in \h^* \, |\, \langle\la ,\al^\vee_i\rangle \in \Z\}$.
A weight $\La \in \h^*$ is dominant integral if $\langle\La, \al_i^\vee\rangle$ are 
non-negative integers for all $i$. 

Fix $\rho\in\h^*$ such that $\langle\rho,\al_i^\vee\rangle=1$,
$i=1,\dots,r$. We have $(\rho,\al_i)= (\al_i,\al_i)/2$.

The Weyl group $\mathcal W\in\on{End (\h^*)}$ is generated by 
reflections $s_i$, $i=1,\dots,r$, 
\be
s_i(\la)=\la-\langle\la,\al_i^\vee\rangle\al_i, \qquad \la\in\h^*.
\ee
We use the notation
\bea\label{shifted}
w\cdot\la=w(\la+\rho)-\rho,\qquad w\in \mathcal W,\;\la\in\h^*,
\eea
for the shifted action of the Weyl group.

The Kac-Moody algebra $\g^t=g(A^t)$ 
corresponding to the transposed Cartan matrix $A^t$ is called
Langlands dual to $\g$.

\subsection{The definition of master functions and critical points} 
Let $\bs\La=(\La_i)_{i=1}^n$, $\La_i\in\mathcal P$;
$\bs z=(z_i)_{i=1}^n\in\C^n$,  $\bs
l=(l_1,\dots,l_r)\in\Z^r_{\geq 0}$,
$\bs  t = ( t_j^{(i)},  j = 1, \dots , l_i)_{i = 1}^r$. 
We call $\La_i$ the {\it weight} at a 
point $z_i$;
$t_j^{(i)}$ a {\it variable of color} $i$.
Define 
\be
\La_\infty\ =\ \sum_{i=1}^n \La_i\  - \ \sum_{i=1}^r
l_i\al_i\ \in\ \mathcal P\ 
\ee
and $\bar {\bs \La}=(\La_1,\dots,\La_n,\La_\infty)$.

The {\it master function} $\Phi_{\g(A)}(\bs t; \bs z; \bar{\bs \Lambda})$
is defined by 
\bean\label{master}
\lefteqn{\Phi_{\g(A)}(\bs t; \bs z; \bar{\bs \Lambda}) = }
\\
&& \prod_{i=1}^r\prod_{j=1}^{l_i}\prod_{s=1}^n
  (t_j^{(i)}-z_s)^{-(\La_s,\al_i)}
  \prod_{i=1}^r\prod_{1\leq j<s\leq l_i} (t_j^{(i)}-t_s^{(i)})^{(\al_i,\al_i)}
  \prod_{1\leq i<j\leq r}\prod_{s=1}^{l_i}\prod_{k=1}^{l_{j}}(t_s^{(i)}-t_k^{(j)})^{(\al_i,\al_j)} \ ,
\notag
\eean
see \cite{SV} .
The function $\Phi$ is a rational function of variables 
$\bs t$ depending on  parameters $\bs z$, $\bs \La$. It is symmetric with
respect to permutations of variables of the same color.


A point $\bs t$ with complex coordinates is called a {\it critical point} of the master function $\Phi$
if the following system of algebraic equations is satisfied
\bean\label{Bethe eqn}
\sum_{s=1}^n \frac{-(\Lambda_s, \alpha_i)}{t_j^{(i)}-z_s}\ +\
\sum_{s,\ s\neq i}\sum_{k=1}^{l_s} \frac{(\alpha_s, \alpha_i)}{ t_j^{(i)} -t_k^{(s)}}\ +\
\sum_{s,\ s\neq j}\frac {(\alpha_i, \alpha_i)}{ t_j^{(i)} -t_s^{(i)}}\
=\ 0\ 
\eean
where $i = 1, \dots , r$, $j = 1, \dots , l_i$.
In other words, a point $\bs t$ is a critical point if
\be
\left(\Phi^{-1}\frac{\partial \Phi }{\partial t_j^{(i)}}\right)(\bs t)=0, \qquad
i=1,\dots,r,\; j=1,\dots l_i.
\ee

The product of symmetric groups
$S_{\bs l}=S_{l_1}\times \dots \times S_{l_r}$ acts on the critical set of the master function
permuting the coordinates of the same color. All orbits have the same cardinality
$l_1! \cdots l_r!$\ .

In the Bethe Ansatz method applied to the Gaudin models \cite{G}
the system of equations (\ref{Bethe eqn}) is called the  Bethe
Ansatz equations. On  relations of critical points  to the  Bethe Ansatz method of
the Gaudin models see \cite{RV}. See also \cite{Fa, Sk1, Sk2, FFR, Fr, MV}.

\subsection{The case of isolated critical points}
In this section we give a sufficient condition for the set of
critical points to be finite and state a conjecture about its cardinality.

We say that the  set of weights $\bar {\bs\La}$ is {\it separating} if
\be
(2\La_\infty\ +\ 2 \rho\ + \ \sum_{i=1}^r c_i\al_i,\ \sum_{i=1}^r c_i\al_i)\
\neq 0,
\ee
for all sets of integers $(c_i)_{i=1}^r$ such that 
$0\leq c_i\leq l_i$, $\sum_i c_i\neq 0$.

For example, if the scalar product is non-negative on the root
lattice and $\La_\infty$ is dominant integral, then $\bar{\bs\La}$ is separating.

The following lemma is a generalization of Theorem 6 in \cite{ScV}.

\begin{lem} If $\bar{\bs\La}$ is separating, then the set of critical
  points is finite. 
\end{lem}

\begin{proof} (Cf. proof of Theorem 6 in \cite{ScV}.)
If the algebraic set of critical points is infinite, then it is
unbounded. Suppose we have a sequence of critical points which is
unbounded. Without loss of generality, we assume that $t_{j}^{(i)}$ tends
to infinity for $i = 1, \dots , r, \ j = 1, \dots , c_i$, and remains bounded for
all other values of $i, j$.

Take the equation \Ref{Bethe eqn} corresponding to a variable
$t_j^{(i)}$ and  multiply it by
$t_j^{(i)}$. Then add the resulting equations corresponding to $i=1,\dots,r$,
$j=1,\dots,c_i$, and pass to the limit along our sequence of critical points.
Then the resulting equation is
\be
(2\La_\infty\ +\ 2\rho\ +\ \sum_{i=1}^r c_i\al_i,\ \sum_{i=1}^r c_i\al_i)\ =\ 0\ .
\ee
This equation  contradicts to our assumption.
\end{proof}

\begin{conj}\label{main conj}
If all components of $\bar{\bs\La}$ are dominant integral weights, 
then for generic $z_1, \dots , z_n$   the number
of $S_{\bs l}$-orbits of critical points of the master function 
$\Phi({\bs{t}; \bs{z}; \bar{\bs{\La}}})$ is equal to 
the multiplicity of the irreducible $\g(A)$-module with highest weight
$\Lambda_\infty$ in the tensor product of irreducible $\g(A)$-modules 
with highest weights $\Lambda_i$, \  $i = 1, \dots , n$. Moreover, all critical points 
are  non-degenerate. \end{conj}

Conjecture \ref{main conj} is proved in \cite{ScV} for 
$\g=sl_2$. In this paper we prove for $\g = sl_n$ 
that the number of $S_{\bs l}$-orbits of critical points
is not greater than the above multiplicity, thus relating the number of critical
orbits and multiplicities of
irreducible representations in tensor products.

Conjecture \ref{main conj} is related to the conjecture on completeness of the Bethe
Ansatz for Gaudin models, see \cite{RV}. In a Gaudin model
to every orbit of critical points one assigns 
an eigenvector in the space of states of a family of commuting linear operators
called  Hamiltonians. The Bethe Ansatz conjecture predicts
that the constructed  eigenvectors span a basis in the space of states.
The dimension of the space of states is equal to the above multiplicity.
Therefore, if the Bethe Ansatz conjecture were true, 
then the number of orbits of critical points would  be not less than the above
multiplicity. The Bethe Ansatz conjecture is proved in \cite{RV} for $\g=sl_2$.

\subsection{On limits of critical points} In this section we formulate auxiliary
results which we use later. 

It is useful to consider
functions more general than master functions,
\be
\Phi (\bs t; \bs z; \bs \mu; \bs \nu ) \ =\ 
\prod_{s=1}^n \prod_{i=1}^{l}
 (t_i - z_s)^{\mu_{i,s}}\
 \prod_{1\leq i < j \leq l} (t_i - t_j)^{\nu_{i,j}} \ .
\ee
For $j > i$, set $\nu_{j, i} = \nu_{i,j}$.
Say that $\bs t = (t_1, \dots , t_l)$ is a critical point of $\Phi$, if
\be
\sum_{s=1}^n \frac{\mu_{i,s}}{t_i - z_s}\ +\
\sum_{j,\ j\neq i} \frac{\nu_{i,j}}{ t_i - t_j}\ = \ 0\ 
\ee
for $i = 1, \dots , l$.  

For any subset $I\subset \{1, \dots , l\}$, denote $\bar I =
\{1, \dots , l\} - I$ its complement. 
Say that the pair $\bs \mu, \bs \nu$ is {\it separating}, if
\be
\sum _{i\in I} \sum _{j\in \bar I} \nu_{i,j}
\ +\ 
\sum _{i, j \in I, \ i<j} \nu_{i,j}
\ +\
\sum _{i\in I} \sum _{s=1}^n  \mu_{i,s}\ \neq \ 0
\ee
for any non-empty subset $I\subset \{1, \dots , l\}$.

\begin{lem}\label{separ}
If the pair $\bs \mu, \bs \nu$ is separating, then the set of critical 
points of $\Phi$ is finite.
\hfill $\square$
\end{lem}

Assume that $\bs t_k = (t_{1,k}, \dots , t_{l,k}), \ k=1, 2, \dots ,$ 
is a sequence of critical points of $\Phi$. Assume that for every $i$,
the sequence $\{ t_{i,k}\}$  has a limit, finite or infinite, as $k$ tends to infinity.
We show that, if we ignore the coordinates whose limit belong to $\{z_1, \dots , z_n, \infty\}$,
then the limits of the remaining coordinates form a critical point of a suitable function
$\tilde \Phi$ defined below.

For $w \in \C \cup \infty$, denote $I(w)$ the subset of $\{1, \dots , l\}$ consisting of all $i$
such that $\lim_{k\to \infty} t_{i,k} = w$.
Let $t_{1,\infty}, \dots t_{\tilde l,\infty}$ be all pairwise distinct limiting points lying in
$\C - \{z_1, \dots , z_n\}$.

Define numbers $\tilde \mu_{\tilde \imath, s}, \ 
\tilde \imath = 1, \dots \tilde l, \ s = 1, \dots , n$, by
\be
\tilde \mu_{\tilde \imath, s} \ =\ \sum _{i \in I(w_{\tilde \imath, \infty})} \mu_{i,s}
\ + \
\sum _{i \in I(w_{\tilde \imath, \infty})} \sum _{j \in I(z_s)} \nu_{i,j}\ .
\ee
Define numbers $\tilde \nu_{\tilde \imath, \tilde \jmath},\  \tilde \imath, \tilde \jmath \in 
\{1, \dots, \tilde l\}, \tilde \imath\neq \tilde \jmath$, by
\be
\tilde \nu_{\tilde \imath, \tilde \jmath} \ =\ \sum _{i \in I(w_{\tilde \imath, \infty})} 
\sum _{j \in I(w_{\tilde \jmath, \infty})} \nu_{i,j }\ .
\ee
Define a function
\be
\tilde \Phi (\tilde{\bs t}; \bs z; \tilde{\bs \mu}; \tilde{\bs \nu} ) \ =\ 
\prod_{s=1}^n \prod_{\tilde \imath = 1}^{\tilde l}
 (\tilde t_{\tilde \imath} - z_s)^{\tilde \mu_{\tilde \imath , s}}\
 \prod_{1\leq \tilde \imath < \tilde \jmath \leq \tilde l} (\tilde t_{\tilde \imath} - \tilde t_{\tilde \jmath})
^{\tilde \nu_{\tilde \imath, \tilde \jmath}} \ .
\ee
\begin{lem}
The point $(t_{1, \infty}, \dots , t_{\tilde l, \infty})$ is a critical point
of the function $\tilde \Phi$.
\hfill $\square$
\end{lem}
\begin{lem}\label{separ-separ}
Assume that the pair $\bs \mu^0 = \{ \mu_{i,s}, \ i \in \{ 1, \dots , l\} - I(\infty),
s=1, \dots , n \}$ and $\bs \nu^0 = \{ \nu_{i,j}, \ i, j \in \{ 1, \dots , l\} - I(\infty),
 i\neq j \}$ is separating. Then the pair
$\tilde {\bs \mu} = \{ \tilde \mu_{\tilde \imath, s}, \ \tilde \imath  = 1, \dots , \tilde l,
s=1, \dots , n \}$ and $\tilde {\bs \nu} = \{ \tilde \nu_{\tilde \imath,\tilde \jmath}, \ \tilde \imath, 
\tilde \jmath \in \{ 1, \dots , \tilde l\},
\tilde  \imath \neq \tilde \jmath \}$ is separating. 
\hfill $\square$
\end{lem}

\section{Populations of critical points}

\subsection{Remarks on Fuchsian equations}

Consider a differential equation  for a function $u(x)$
\bean\label{fuchs eqn}
u^{(k)}\ +\ p_1 u^{(n-1)}\ +\ \dots\ +\ p_k u\ =\ 0\ ,
\eean
where $p_i = p_i(x)$ are rational functions. A point $z\in\C$ is
called an {\it ordinary point} if all $p_i(x)$ are  holomorphic at $z$.
A non-ordinary point is called {\it singular}. 

A singular point $z \in \C$ is called {\it regular} if the order of the pole of
$p_i$ at $z$ is at most $i$. 

Equation \Ref{fuchs eqn} has an ordinary (resp., singular,
regular singular) 
point at infinity if after the change of variable $x = 1/\xi$ the
point $\xi = 0$ is ordinary (resp., singular,
regular singular).

A differential equation with only regular singular points is called
{\it Fuchsian}. 

Let $f_1, \dots , f_k$ be linearly independent
polynomials. There is a unique (up to multiplication by a function) 
linear differential equation of order
$k$ with solutions $f_1, \dots , f_k$,
\be
W(u, f_1, \dots , f_k) \ / \ W(f_1, \dots , f_k)\ =\ 0\ ,
\ee
where  
\be
W(g_1, \dots , g_s)\ =\ \det(g_i^{(j-1)})_{i, j=1}^s
\ee
is the {\it Wronskian} of functions $g_1, \dots , g_s$. 
This equation is Fuchsian. 

Consider a Fuchsian equation  \Ref{fuchs eqn} and write in
a neighborhood of a point $z\in\C$
\be
p_i\ =\ \sum_{s=0}^\infty \ p_{is}\ (x-z)^{s-i}, \qquad i = 1, \dots , k.
\ee
If a function
\be
u\ =\ (x-z)^\la \ (1\ +\ \sum_{s=1}^\infty a_s\ (x-z)^s)
\ee
is a solution of equation \Ref{fuchs eqn}, 
then $\la$ is a root of the {\it indicial equation} at
the point $z$
\bea\label{indicial eqn}
\la (\la-1) \dots (\la-k+1)\ +\ p_{10} \la (\la-1) \dots (\la-k+2)\ +
\ \dots\ +\ p_{k0}\ =\ 0 .
\eea
The roots of the indicial equation at a point $z$ are called {\it
  exponents} at $z$ of the Fuchsian equation. If $V$ is the space of solutions of 
the Fuchsian equation,
then the roots of the indicial equation are called the exponents of $V$ at $z$.
 

\subsection{ Polynomials representing critical points }\label{PLCP}
Let $\bs t=(t_j^{(i)})$ be a critical point of a master function 
$\Phi=\Phi({\bs{t};\bs{z}; \bar{\bs{\La}}})$.
Introduce polynomials $\bs y=( y_1(x),\dots, y_r(x))$,
\be
y_i(x)\ =\ \prod_{j=1}^{l_i}(x-t_j^{(i)}).
\ee
The $r$-tuple $\bs y$ determines the $S_{\bs l}$-orbit of the critical point $\bs t$.
We say that the $r$-tuple of polynomials $\bs y$ {\it represents a 
critical point} of the master function $\Phi$. Usually we do not distinguish between
the critical points of the same $S_{\bs l}$-orbit and identify a critical point with
the $r$-tuple of polynomials $\bs y$ representing the point.

We consider the $r$-tuple $\bs y$ up to multiplication 
of each coordinate by a non-zero number,  since we are interested only in the roots of polynomials
$y_1, \dots , y_r$. Thus the $r$-tuple defines a point in the direct product 
$\PCr$ of $r$ copies of the projective space associated with the vector 
space of polynomials of $x$.

Introduce polynomials
\bean\label{T}
{T}_i(x)=\prod_{s=1}^n(x-z_s)^{\langle \La_s, \al_i^\vee\rangle}, \qquad i = 1, \dots , r .
\eean

We say that a given $r$-tuple of polynomials $\bs y \in \PCr$ is {\it generic with respect 
to integral dominant weights  $\Lambda_1, \dots , \Lambda_n$ of the 
Kac-Moody algebra $\g(A)$ and points}
$z_1, \dots z_n$ if
\begin{itemize} 
\item
each polynomial $y_i(x)$ has no multiple roots;
\item
all roots of $y_i(x)$ are different from roots of the polynomial $T_i$;
\item
any two polynomials $y_i(x)$, $y_j(x)$ have no common roots if $i\neq j$ and
$a_{ij}\neq 0$.
\end{itemize}

If $\bs y$ represents  a critical point of $\Phi$,  then $\bs y$ is generic.

Now we reformulate the property of $\bs y$ to represent a critical point. 

Write $f' = \partial f/\partial x$ and $\ln'(f) = f'/f$. 
Let polynomials $F_i$, $G_i, \ i = 1, \dots , r$, be given by
\be 
F_i\ = \ \prod_{s=1}^n (x-z_s)\ \prod_{j,\ a_{ij}\neq 0} y_j, \qquad 
G_i\ =\ F_i \ \ln'\left(T_i \prod_{j,\ j\neq i} y_j^{-\langle\al_j,\al_i^\vee\rangle}\right).
\ee

\begin{lem}\label{deg 2 lem} 
A generic $r$-tuple $\bs y$ represents a critical point if and only if 
for every  $i = 1, \dots , r$ the polynomial
$F_i \ y_i''\  - \ G_i\  y_i'$ is  divisible by the polynomial $y_i$.  In other
words, a generic $r$-tuple $\bs y$ represents a critical point if and only if 
for every  $i = 1, \dots , r$ there exists a
polynomial $H_i$, such that deg $H_i \ \leq $ deg $F_i\ -\ 2$, and
the polynomial $y_i$ is a solution of the differential equation
\bean\label{deg 2 eqn}
F_i\ u''\ -\ G_i\ u'\ +\ H_i\  u\ =\ 0\ . 
\eean
\end{lem}
\begin{proof} The lemma is a direct corollary of a classical result
  of Heine-Stieltjes, see \cite{S}, Section 6.8. We sketch the proof.

Assume that there exist such polynomials $H_1, \dots , H_r$. 
Substitute $x=t_j^{(i)}$ into $F_i \ y_i''\ -\ G_i\  y_i'\ +\ H_i \ y_i\ =\ 0$. We get 
$y''(t_j^{(i)}) / y_i'(t_j^{(i)})\ $ $ =\ $ $G_i(t_j^{(i)}) / F_i (t_j^{(i)})$. This is exactly equation
\Ref{Bethe eqn} multiplied by $2/(\al_i,\al_i)$, since for
$f = \prod_s(x-a_s)$ we have
\be
\frac{f''}{f'}(a_s)\ =\ \sum_{k, \ k\neq s}\ \frac{2}{a_s-a_k}\ .
\ee
This means that the roots of polynomials $y_1, \dots , y_r$ form a critical point. 
This argument is reversible.
\end{proof}

Let $\bs y$ represent a critical point. Then  equation \Ref{deg 2
  eqn} is  Fuchsian. The singular points and exponents of that equation are 
\bea
\begin{array}{ll}
x=z_s:\qquad & \{0,\langle\La_s,\al_i^\vee\rangle+1\},\\
x=t_j^{(k)},\; k\neq i: \qquad & \{0,-\langle\al_k,\al_i^\vee\rangle+1\},\\
x=\infty : \qquad & \{-l_i, -\sum_s\langle\La_s,\al_i^\vee\rangle+\sum_{k,\ k\neq i}
l_k\langle\al_k,\al_i^\vee\rangle-1\}. 
\end{array}
\eea

\begin{lem} \label{wronskian-critical lem}
Let $\bs y$ be generic and let $\langle\La_s,\al_i^\vee\rangle$ be non-negative
integers for some $i$ and all
$s=1, \dots , n$. Then $F_i\ y_i''\ -\ G_i\ y_i'$ is divisible by
$y_i$ if and only if there exists a polynomial $\tilde {y}_i(x)$ such
that the  Wronskian $W(y_i, \tilde{y}_i)$ is given by
\bean\label{wronskian-critical eqn}
W(y_i,\tilde{y}_i)\ = \ T_i \ \prod_{j,\ j\neq i}y_j^{-\langle\al_j,\al_i^\vee\rangle}\ .
\eean
\end{lem}


\begin{proof}
If $H_i\ = \ (F_i\ y_i''\ -\ G_i\ y_i')/ y_i$ is a polynomial, then
by \cite{ScV}, Lemma 7, all solutions of $F_i\  u'' \ - \  G_i\ u'\ $ $ +\ H_i\ u\ =\ 0$ are
polynomials. Then $\tilde{y}_i$ is any second linearly independent solution
multiplied by a suitable constant. This proves the ``only if'' part of
the lemma. 

Let a polynomial $\tilde{y}_i$ exist.  The polynomial $y_i$ is a solution of
the equation $W(u, y_i, \tilde{y}_i) = 0$. After multiplying this equation by
$F_i /W(y_i,\tilde{y}_i)$ we get $F_i u'' - G_i u' + H_i  u = 0$, where
\be
H_i(x)\ =\
F_i\;\
\frac{y_i'\ \tilde y_i''\ -\ y_i''\ \tilde
  y_i'}{W(y_i, \tilde y_i)}\ =\
F_i\;\ \frac{\ln'( W(y_i, \tilde
  {y}_i))\ y_i'-y_i''}{y_i} \ .
\ee
It is clear that 
$F_i \ln'( W(y_i, \tilde y_i) )$ is a polynomial. Therefore, poles of
$H_i$ are common zeros of $y_i$ and $W(y_i, \tilde y_i)$. 
Equation \Ref{wronskian-critical eqn} implies that the 
polynomials $y_i$ and $W(y_i, \tilde y_i)$ do not have common zeros
since $\bs y$ is assumed to be generic.  Therefore $H_i$ is a
polynomial. 
\end{proof}

\begin{cor}\label{fertile cor}
Let the weights $\Lambda_1, \dots , \Lambda_n$ be dominant integral. Then
a generic $r$-tuple $\bs y$ represents a critical point
if and only if for every $i = 1, \dots , r$
there is a polynomial $\tilde y_i$ satisfying 
(\ref{wronskian-critical eqn}).
\end{cor}

\begin{lem}\label{tilde y}
Let $y_1, \dots, y_r,\ T_i$ be given and let $\tilde y_i$ satisfy equation
\Ref{wronskian-critical eqn}. Then, up to multiplication by a non-zero number,
the function $\tilde y_i$ has the form
\bean\label{tilde}
\tilde y_i(x)\ =\ c_1\ y_i (x) \ \int \ T_i(x)\ \prod_{j=1}^r 
y_j^{-\langle\al_j,\al_i^\vee\rangle}
\ dx \ +\ c_2\ y_i(x)\ ,
\eean
where $c_1, c_2$ are complex  numbers.
\hfill $\square$
\end{lem}
Notice  that formula \Ref{tilde} gives all solutions of the  differential equation
\be
F_i\ u''\ -\ G_i\ u'\ +\ H_i\  u\ =\ 0\ .
\ee

Lemma \ref{tilde y} shows that the $r$-tuples 
\bean\label{simple}
\bs y^{(i)}\ =\ (y_1, \dots , \tilde y_i,\dots, y_r) \ 
\qquad \in \ \qquad \PCr \ ,
\eean
where $\tilde y_i$ is given  by \Ref{tilde}, form a one-parameter family.
The parameter space of the family
is identified with the projective line $\bs P^1$ with projective coordinates 
$(c_1 : c_2)$.
We have a map
\bean\label{map}
Y_{\bs y, i}\ : \ \bs P^1 \ \to \PCr\ , 
\eean
which sends a point $c = (c_1 : c_2)$ to the corresponding $r$-tuple $\bs y^{(i)}$.

\begin{lem}\label{generic lem}
If $\bs y$ is generic, then almost all $r$-tuples $\bs y ^{(i)}$ are generic.
The exceptions form a finite set in $\bs P^1$.
\hfill $\square $
\end{lem}

\subsection{Fertile $r$-tuples}\label{FER}

Let $\La_1, \dots, \La_n$  be dominant integral weights, $z_1, \dots , z_n$
complex numbers. Let $\bs y = (y_1, \dots , y_r)\ \in\ \PCr$ and
let $l_i $ be the degree of the polynomial $y_i$.
The weight 
\be
\La_\infty \ =\  \sum_{s=1}^n\ \La_s\  -\  \sum _{i=1}^r\  l_i\  \al_i
\ee 
is called {\it the weight at infinity} of the $r$-tuple $ {\bs y}$ with respect to
$\La_1, \dots, \La_n$  and $z_1, \dots , z_n$.

An $r$-tuple $\bs y$  is called {\it fertile} 
with respect to  dominant integral weights $\La_1, \dots , \La_n$ and  
complex numbers $z_1, \dots , z_n$,  if for every $i$ there
exists a polynomial $\tilde y_i$ satisfying equation
\Ref{wronskian-critical eqn}. If $\bs y$ is fertile, then the $r$-tuples
$\bs y^{(i)}$ given by \Ref{simple} are called {\it immediate descendants} of $\bs y$
in the $i$-th direction.

A generic $r$-tuple $\bs y$ represents a critical point of a master 
function associated to dominant integral weights $\Lambda_1, \dots  , \Lambda_n$ 
and points $z_1, \dots , z_n$ if and only if it is fertile, see
Corollary \ref{fertile cor}.

\begin{lem}\label{limit fertile lem}
Assume that a sequence $\bs y_k$,\ $k = 1, 2, \dots$ ,\ of fertile $r$-tuples of 
polynomials has a limit $\bs y_\infty$ in $\bs P (\C[x])^r$ as $k $ tends to infinity.
\begin{enumerate}
\item[$\bullet$] Then the limiting $r$-tuple $\bs y_\infty$ is fertile.
\item[$\bullet$] Let $i \in \{ 1, \dots  , r\}$.
Let $\bs y^{(i)}_\infty$ be an immediate descendant of $\bs y_\infty$ 
in the $i$-th direction. Then for any $k$, 
there exists an immediate descendant
$\bs y_k^{(i)}$ of $\bs y_k$ 
such that 
$\bs y^{(i)}_\infty$ is the limit of $\bs y_k^{(i)}$ as $k$ tends to infinity.
\end{enumerate}
\end{lem}

\begin{proof}
Let $\bs y_k = (y_{k,1}, \dots , y_{k,r})$. For every $k$, including $k=\infty$, 
consider the differential equation 
\bean\label{local} 
F_{k,i}\ u''\ -\ G_{k,i}\ u'\ +\ H_{k,i}\  u\ =\ 0 ,
\eean
where $F_{k,i}$, $G_{k,i}$ are as in \Ref{deg 2 eqn} with $y_j$ replaced by $y_{k,j}$
and $H_{j,i}\  =\ (G_{k,i} y_{k,i}'- F_{k,i} y_{k,i}'') / y_{k,i}$.
For every $k < \infty$, all solutions of that equation are polynomials, and the polynomial
$y_{k,i}$ is one of solutions. 

The lemma would follow if we proved that all solutions of equation \Ref{local} for $k=\infty$
were polynomials.

Since the sequence of $\bs y_k$ has a limit,  there is a point $z \in\C$, 
such that $x=z$ is an ordinary point of equation \Ref{local} for all $k$.
 Fix $a, b \in \C$, and let $\tilde y_{k,i}$ be the solution of 
\Ref{local} with the initial condition 
$\tilde y_{k,i}(z) = a$,\  $\tilde y_{k,i}'(z) = b$.

By the standard theorem on continuous dependence of solutions on the coefficients of the equation, 
the function $\tilde y_{\infty, i}$ is the limit of functions
$\tilde y_{k,i}$ as $k$ tends to infinity.

The function $\tilde y_{\infty,i}$ is univalued and regular,  since $\tilde y_{k,i}$ 
is a polynomial for every finite $k$. Hence the function $\tilde y_{\infty,i}$
is a polynomial. Thus all solutions of equation \Ref{local} for $k = \infty$ 
are polynomials. This implies the lemma.
\end{proof}

Let $\bs y$ represent a critical point of a master function
 $\Phi(\bs t;\bs z; \bs \La, \La_\infty)$. 
Let ${\bs y}^{(i)}$ $=$ $(y_1,$ $ \dots , $ $ \tilde y_i,$ $ \dots, y_r)$, 
be an immediate descendant of $\bs y$ and let
$\La^{(i)}_\infty$ be its weight at infinity with respect to 
$\La_1, \dots, \La_n$  and $z_1, \dots , z_n$.

The following observation is crucial in  this paper.

\begin{theorem}\label{duplicate thm}
If ${\bs y}^{(i)}$ is generic, then ${\bs y}^{(i)}$ represents a critical
point of the master function $\Phi(\bs t;\bs z; \bs \La,  \La^{(i)}_\infty)$.
\end{theorem}

\begin{proof}
Denote  $\tilde t_s^{(i)}$ the   roots of $\tilde y_i$. 
For any $j$, such that $j \neq i$ and $a_{ij}\neq 0$, choose a root $t^{(j)}_k$ of the polynomial
$y_j$. We have $W(y_i, \tilde y_i)(t^{(j)}_k) = 0$ by (\ref{wronskian-critical eqn} ).
Hence
\be
\sum_{s}\frac{1}{t_k^{(j)}-t_s^{(i)}}\ =\ \sum_{s}\frac{1}{t_k^{(j)}-\tilde
  t_s^{(i)}}\ .
\ee
This implies that the roots of $\bs y^{(i)}$ satisfy the equation of system
(\ref{Bethe eqn}) corresponding to the coordinate $t^{(j)}_k$.

The roots of $\bs y^{(i)}$ satisfy the equations of system 
(\ref{Bethe eqn}) corresponding to coordinates $\tilde t^{(i)}_s$
according to Lemma \ref{wronskian-critical lem}. 
\end{proof}

Thus, starting with an $r$-tuple $\bs y $
representing a critical point of a master function
$\Phi(\bs t;\bs z; \bs \La, \La_\infty)$ and an index 
$i \in \{1, \dots , r\}$, we  construct in  \Ref{map} a  family 
$Y_{\bs y, i} : \bs P^1 \to \PCr$ of fertile $r$-tuples.
For almost all $c \in \bs P^1$ (with only finitely many exceptions), 
the $r$-tuple $Y_{\bs y, i}(c)$  represents a critical point of a master function
associated with integral dominant weights $\La_1, \dots , \La_n$ and points
$z_1, \dots , z_n$.

We call this construction the {\it simple reproduction procedure in the $i$-th direction}.

\subsection{General reproduction procedure}\label{general procedure}

Assume that the weights $\La_1, \dots , \La_n$ are dominant integral and
an $r$-tuple $\bs y^0 \in \PCr$ represents a critical point of a master function
$\Phi(\bs t; \bs z; \bs \La, \La_\infty)$. 

Let $\bs i = (i_1, i_r, \dots , i_k), \ 1 \leq i_j \leq r,$ be a sequence of natural numbers.
We define a $k$-parameter family of fertile $r$-tuples 
\bea\label{general map}
Y_{\bs y^0, \bs i} \ :\ (\bs P^1)^k \ \to \PCr
\eea
by induction on $k$,  starting at $\bs y^0$ and successively applying 
the simple reproduction procedure in directions
$i_1, \dots , i_k$.

More precisely, for $k=1$, it is the family 
$Y_{\bs y^0, i_1}  : \bs P^1 \ \to \PCr$ defined by \Ref{map}.
If $k > 1$, \ $\bs i' = ( i_1, i_2, \dots , i_{k-1})$,\ 
$\bs c = (c^1 , \dots , c^k) \in (\bs P^1)^k$, and
$\bs c' = (c^1 , \dots , c^{k-1}) \in (\bs P^1)^{k-1}$,
then we set
\be
Y_{\bs y^0, \bs i} ( \bs c) \ = \
Y_{Y_{\bs y^0, \bs i'}(\bs c'), i_k} (c^k)\ .
\ee
The image $P_{\bs y^0, \bs i} \ \subset \ \PCr$ of the map 
$Y_{\bs y^0, \bs i}$ is called the {\it population in the direction of} 
$\bs i$ {\it originated at} $\bs y^0$.
The set $P_{\bs y^0, \bs i}$ is an irreducible algebraic variety.

It is easy to see that if $\bs i' = (i'_1, i'_r, \dots , i'_{k'}), \ 1 \leq i'_j \leq r$,
is a sequence of natural numbers, 
and the sequence $\bs i'$ is contained in the sequence $\bs i$ as  an ordered subset, then
 $P_{\bs y^0, \bs i'}$ is a subset of $P_{\bs y^0, \bs i}$.

The union 
\be 
P_{\bs y^0} \ = \ \cup_{\bs i} \ P_{\bs y^0, \bs i}\ \subset \PCr\ ,
\ee
where the summation is over all sequences $\bs i$, 
is called {\it the population of critical points associated} with the Kac-Moody algebra $\g$, 
weights $\La_1, \dots , \La_n$, points $z_1, \dots , z_n$,  and {\it originated} at $y^0$.

\begin{lem}\label{open}
For a given $\bs i = (i_1, \dots , i_k)$, almost all $r$-tuples 
$Y_{\bs y^0, \bs i} ( \bs c)$ represent  critical points of master functions associated to
weights $\La_1, \dots , \La_n$, and points $z_1, \dots , z_n$.
 Exceptional values of $c \in (\bs P^1)^k$ 
are contained in a proper algebraic subset.
\hfill $\square$
\end{lem}

\begin{lem}\label{intersect}
If two populations intersect, then they coincide.
\hfill $\square$
\end{lem}

\subsection{Populations and flag varieties}\label{Conjec}

${\ }$

{\bf Example.} 
Consider the population of critical points associated to
$\g = sl_3$ and $n=0$ and originated at $\bs y^0 = (1, 1)$. The pair $(1, 1)$
represents the critical point of the function with no variables. 
This population consists of pairs of non-zero 
polynomials $\bs y = (y_1, y_2)$, where 
\bean\label{Example}
y_i\ = \ a_{2,i}\ x^2\ + \ a_{1,i}\ x\ + \ a_{0,i}\ , \qquad i=1, 2 \ ,
\eean
and
\be
a_{1,1} \ a_{1,2}\ =\ 2\ a_{0,1} \ a_{2,2}\ + 2\ a_{2,1} \ a_{0,2}\ .
\ee
 For any  pair $\bs y = (y_1, y_2)$, 
if $y_1, y_2$ do not have multiple
roots and do not have common roots, then the roots of the 
polynomials $y_1, y_2$ 
form a critical point of the function
\be
\Phi \ = \ \prod_{ 1\leq i < j \leq l_1 } (t^{(1)}_i - t^{(1)}_j )^2\
\prod_{ 1\leq i < j \leq l_2 } (t^{(2)}_i - t^{(2)}_j )^2\
\prod_{ i=1}^{ l_1 } \prod_{ j=1}^{ l_2 } (t^{(1)}_i - t^{(2)}_j )^{-1} \, , 
\ee
where $l_1 = $ deg $y_1$ and $l_2 = $ deg $y_2$.

In this case equations (\ref{wronskian-critical eqn}) take the form
\bean\label{example}
W(y_1, \tilde{y}_1)\ = \ y_2\ , \qquad W(y_2, \tilde{y}_2)\ = \ y_1\ ,
\eean
and the reproduction procedure works as follows. We start with $\bs
y^0 = (1, 1)$. 
Equations (\ref{example}) have the form 
$W( 1,\ \tilde{y}_1) =  1$, \ $ W( 1,\ \tilde{y}_2) =   1$.
Using, the first of them, we get pairs $\bs y = (x + a, \ 1)$ for all numbers $a$. 
Equations (\ref{example}) now are
$W( x + a, \ \tilde{y}_1) =  1$, \ $ W( 1,\ \tilde{y}_2) =  x + a$.
Using the second equation we get pairs $\bs y = ( x + a, \ x^2 / 2 + ax + b)$ for all  $a, b$.
Equations (\ref{example}) take the form 
$W( x + a, \ \tilde{y}_1) =  x^2/2 + ax + b$, \ $W(  x^2/2 + ax + b, \ \tilde{y}_2) =  x + a$.
Using the first of them  we get 
$\bs y = (x^2/2 + cx + ac-b, \ x^2/2 + ax + b)$ for all  $a, b, c$. 

If we started the procedure using equation $ W( 1,\ \tilde{y}_2) =   1$,
then the constructed pairs would have been of the form $\bs y = (1,\ x + a)$,
$\bs y = (x^2 / 2 + ax + b, \ x + a)$, $\bs y = (x^2/2 + ax + b,\ x^2/2 + cx + ac-b)$.

It is easy to see that the
union of all those pairs is our population, and nothing else can be
constructed 
starting from $\bs y^0 = (1, 1)$. 

\medskip

It is easy to see that the family of pairs (\ref{Example}) (where 
each pair is considered up to multiplication of its coordinates by 
non-zero numbers) is isomorphic as an algebraic variety to the variety of all
full flags in the three dimensional vector space $V$ of the first
coordinates of the pairs. 
Namely, $y_1$ generates a line in $V$ and $y_2$ defines a plane in 
$V$ containing the line generated by $y_1$.

\begin{conj}
Every population, associated to a Kac-Moody algebra $\g$, dominant integral weights
$\Lambda_1, \dots , \Lambda_n$,   points $z_1, \dots , z_n$,
 is a (pro-)algebraic variety isomorphic to 
the full flag variety associated to the Kac-Moody algebra
$\g^t$ which is Langlands dual to $\g$. Moreover, the parts of the family corresponding to 
$r$-tuples of polynomials with fixed degrees are isomorphic to 
Bruhat cells of the flag variety.
\end{conj}

In the example above the possible degrees of polynomials $y_1, y_2$ are
(0,0), (1,0), (0,1), (1,2), (2,1), (2,2). The corresponding parts of the family
are isomorphic to open Bruhat cells of dimensions 0, 1, 1, 2, 2, 3, respectively.

In this paper we prove this conjecture for the Lie algebras with root systems 
of types $A_N, B_N, C_N$.

\subsection{Degrees of polynomials in a population and the Weyl group}

Let $\La_1, \dots, \La_n$  be dominant integral weights, $z_1, \dots , z_n$
complex numbers. Let an $r$-tuple $\bs y$ be a member of a population $P$
associated to dominant integral weights $\La_1, \dots, \La_n$ and 
points $z_1, \dots , z_n$. The $r$-tuple $\bs y$ is fertile by Lemma \ref{limit fertile lem}.
Let an $r$-tuple $ {\bs y}^{(i)} = (y_1, \dots , \tilde y_i, \dots , y_r)$ be 
an immediate descendant of $\bs y$.  
Let  $ \La^{(i)}_\infty$ be the weight at infinity of $ {\bs y}^{(i)}$.

\begin{lem}\label{simple inf lem} 
If the degree of $\tilde y_i$ 
is not equal to the degree of  $y_i$, then  
\be
\Lambda^{(i)}_\infty\ = \ s_i\cdot{\La}_\infty\ ,
\ee
where $s_i\cdot $ is the shifted action of the $i$-th generating
reflection of the Weyl group.
\end{lem}
\begin{proof}
Formula \Ref{wronskian-critical eqn} implies 
\be
\tilde l_i\ +\ l_i\ -\ 1\ =\ \deg \ W(y_i, \tilde y_i)\ =\
\langle \ \sum_s\ \La_s\ -\ \sum_{j,\ j\neq i}\ l_j\al_j,\ \al_i^\vee\ \rangle \ .
\ee
Therefore
\bea
(\ \tilde l_i\ -\ l_i\ )\ \al_i\ &=&\ \langle \ \sum_s
\ \La_s\ -\ \sum_{j}\ l_j\al_j\ +\
\rho,\al_i^\vee\rangle \ \al_i\ 
\\
&=&
 \La_\infty\ +\ \rho\ - \ \ s_i\ (\ \La_\infty\ +\ \rho\ )\
=\ \La_\infty\ - \ s_i\cdot\ \La_\infty\ .
\eea
\end{proof}

\begin{theorem}\label{inf weight thm}
Let $\La_1, \dots, \La_n$  be dominant integral weights.
Let $\bs y^0$ be a critical point of a master function 
$\Phi(\bs t;\bs z; \bs \La, \La_\infty)$. 
Let $P$ be the population
of critical points associated to weights $\La_1, \dots, \La_n$,
points $z_1, \dots , z_n$, and originated at $\bs y^0$.
Then
\begin{enumerate}
\item[$\bullet$]
For any $r$-tuple $\bs y \in P$,  there is an element $w$  of the Weyl group
$\mathcal W$, 
such that the weight at infinity of $\bs y$
is  $w \cdot \Lambda_\infty$.
\item[$\bullet$] 
For any element $w \in  \mathcal W$,  there is an $r$-tuple $\bs y \in P$ whose weight at infinity
is  $w \cdot \Lambda_\infty$.
\end{enumerate}
\hfill $\square$
\end{theorem}

\begin{cor}
If the Weyl group is finite, then every population is a finite dimensional algebraic variety.
\end{cor}

More generally, for a non-negative integer $d$,
let $\C_d[x] \subset \C[x]$ be the vector subspace of polynomials of degree
not greater  than $d$.

\begin{cor}
For every $d$ and  every population $P$, associated to 
a Kac-Moody algebra, dominant integral weights
$\La_1, \dots , \La_n$, points $z_1, \dots , z_n$, 
the intersection $P \cap \bs P (\C_d[x])^r$
is an algebraic variety.
\end{cor}

Theorem \ref{inf weight thm} gives some sufficient conditions for 
a master function not to have critical points.

\begin{cor}
Let $\Phi$ be the master function associated to dominant integral weights 
$\La_1,  \dots , \La_n$ and a weight at infinity $\Lambda_\infty$.
Assume that there is an element $w$ 
of the Weyl group such that $\sum_{s=1}^n \La_s\ -\ w \cdot  \La_\infty$ 
does not belong to the cone $\Z_{\geq 0} \al_1\ \oplus\ \dots\ \oplus \ 
\Z_{\geq 0} \al_r$. Then the master function $\Phi$
does not have critical points.
\end{cor}

\begin{cor}
Let $\Phi$ be the master function associated to dominant integral weights 
$\La_1,  \dots , \La_n$ and a weight at infinity $\Lambda_\infty$.
Assume that  $\La_\infty$ belongs to one of the reflection hyperplanes 
of the shifted action of the Weyl group.
 Then the master function $\Phi$
does not have critical points.
\end{cor}

The next corollary says that under certain conditions on weights there is 
only one population of critical points.

Namely, the $r$-tuple $(1, \dots , 1) \in \PCr $ is the unique $r$-tuple of non-zero polynomials
of degree 0. The weight at infinity of $(1, \dots , 1)$ is
$\Lambda_{\infty, (1, \dots , 1)} = \sum_{s=1}^n \La_s$. Let $O_{(1, \dots , 1)}$ be the orbit of the weight
$\Lambda_{\infty, (1, \dots , 1)}$ with respect to the shifted action of the Weyl group.
Let $P_{(1, \dots , 1)}$ be the population associated  to dominant integral weights 
$\La_1,  \dots , \La_n$, complex numbers $z_1, \dots, z_n$, and originated at $(1, \dots , 1)$.

\begin{cor}
Let $\Phi(\bs t;\bs z; \bs \La, \La_\infty)$  be a master function 
such that $\La_\infty \in O_{(1, \dots , 1)}$. If the master function has
critical points, then they belong to the population  $P_{(1, \dots , 1)}$.
\end{cor}

\subsection{Finiteness of the number of populations}
\begin{theorem}\label{finiteness}
Let $\Phi (\bs t; \bs z; \bs \La, \La_\infty)$ be the master
 function associated to dominant integral
weights $\La_1, \dots \La_n$, a weight at infinity $\La_\infty$, complex numbers
$z_1, \dots , z_n$. Then the number of populations originated at critical points of
the function $\Phi$ is finite, if the orbit of $\La_\infty$ with respect to the 
shifted action of the Weyl group contains a weight $w\cdot \La_\infty$ such that 
the set $\bs \La, \ w\cdot \La_\infty$ is separating.
\end{theorem}
The theorem is a direct corollary of Lemmas \ref{separ} - \ref{separ-separ}.

\section{Vector Spaces of Polynomials and Schubert Cells}

In this section we remind known facts about Schubert calculus.
 
Let $\mc V$ be a complex vector space of dimension $d + 1$  and
\be
\mathcal F=\{ 0\subset F_1\subset F_2\subset\dots\subset
F_{d+1} =\mc V\}, \qquad \dim F_i=i, 
\ee
 a full flag in $\mc V$.  Let $Gr(N+1,\mc V)$ be the Grassmannian
of all $N+1$ dimensional subspaces in $\mc V$.

Let $\bs a =(a_1,\dots,a_{N+1})$, 
$ d - N \geq a_1\geq a_2\geq \dots \geq a_{N+1}\geq 0$,
be a non-increasing sequence of  non-negative integers.
Define the  {\it Schubert cell} $G^0_{\bs a}(\mc F)$
assocaiated to the flag $\mathcal F$ and sequence $\bs a$ as the set
\bea
 \{V\in Gr(N+1,\mc V)\ | &&
\dim (V \cap F_{d- N + i-a_i}) = i ,
\\
&& \dim (V \cap F_{d- N + i-a_i-1}) = i-1, 
\ \text{for}\ 
i = 1, \dots , N+1 \} .
\eea
The closure $G_{\bs a}(\mc F)$ of the Schubert cell is called {\it
the Schubert cycle}. For a fixed flag $F$, the Schubert cells 
form a cell decomposition of the Grassmannian.
The codimension of $G^0_{\bs a}(\mc F(z))$
in the Grassmannian is $|\bs a| = a_1 + \dots + a_{N+1}$. The cell corresponding to
$\bs a = (0, \dots , 0)$ is open in the Grassmannian.

Let $\mc V = \C_d[x]$ be the space of polynomials
of degree not greader than $d$, dim $\mc V = d + 1$. For any  $z \in \C \cup \infty$,
define a full flag in $\C_d[x]$,
\be
\mathcal F(z)\ = \ \{ 0 \subset F_1(z) \subset F_2(z) \subset \dots \subset
F_{d+1}(z) \}\ .
\ee 
For $z \in \C $ and any $i$, let $F_i (z)$ be the subspace of all
polynomials divisible by $( x - z )^{d+1-i}$. For any $i$,
let $F_i(\infty) $ be the subspace of all
polynomials of degree less than $i$.

Thus,  for any sequence $\bs a$ and any $z \in \C \cup \infty$, 
we have a Schubert cell $G^0_{\bs a}(\mc F(z))$ in the Grassmanian
$Gr (N+1, \C_d[x])$ of all $N+1$-dimensional subspaces of $\C_d[x]$.

Let $V \in Gr (N+1, \C_d[x])$.
For any $z \in \C \cup \infty$, let $\bs a(z)$ be such a unique
sequence  that $V$ belongs to the cell $G^0_{\bs a (z)}(\mc F(z))$.
We say that a point $z \in \C \cup \infty$ is {\it a ramification point } for $V$, 
if  $\bs a(z) \neq (0, \dots , 0)$.

This standard terminology is motivated by the fact that each $V \in Gr (N+1, \C_d[x])$ 
defines a rational map $\bs P^1 \to \bs P^N$.

\begin{lem}\label{codim}
For a basis $u_1, \dots , u_{N+1}$  in $V$, let
\be
W(u_1, \dots , u_{N+1}) \ =\ c \ \prod_{s=1}^n(x \ - \ z_s)^{m_s}\ , \qquad c\neq 0\ ,
\ee
be the Wronskian. Then
\begin{enumerate}
\item[$\bullet$] The ramification points for $V$ are the points $z_1, \dots , z_n$ and possibly
$\infty$.
\item[$\bullet$] We have $| \bs a (z_s)| = m_s$ for every $s$.
\item[$\bullet$] We have $| \bs a (\infty)| = (N+1)(d - N) - \sum_{s=1}^n m_s $.
\end{enumerate}
\hfill $\square$
\end{lem}

\begin{cor}\label{sum of cod} (Pl\"ucker formula)\newline
We have
\bean\label{ram cond}
\sum _{s=1}^n | \bs a (z_s)| \ + \ | \bs a (\infty)|\ =\ 
\dim \ Gr ( N+1, \C_d[x]) \ .
\eean
\end{cor}

A point $z \in \C$ is called {\it a base point} for $V$ if $u(z) = 0$ for every $u \in V$.

Assume that ramification conditions are fixed at $z_1, \dots , z_n, \infty$ so that 
\Ref{ram cond} is satisfied.

{\bf Enumerating  Problem.} What is the number of $V$ with those ramification properties?

By the Brill-Noether Theorem in \cite{EH}, the number of such $V$ is always finite. The estimate 
from above is given by the Schubert calculus. The intersection of Schubert cycles
is a union of Schubert cycles, 
\bean\label{cycle mult}
G_{\bs a}\ \cdot\ G_{\bs b} \ = \ \sum_{\bs c}\  N_{\bs a, \bs b, \bs c}\ G_{\bs c} \ ,
\eean
where the cycles are considered up to algebraic equivalence  \cite{Fu}. The structure constants of
this multiplication are determined by the Littlewood-Richardson rule, which has 
a representation theoretic sense.

Namely, for a non-increasing sequence $\bs a =(a_1,\dots,a_{N+1})$, 
$  a_1\geq a_2\geq \dots,\geq a_{N+1}\geq 0$,
of  non-negative integers, denote
$\tilde L_{\bs a}$ the finite dimensional irreducible $gl_{N+1}$-module with highest weight $\bs a$.
The tensor product of irreducible modules is a direct sum of irreducible modules,
\bean\label{rep mult}
\tilde L_{\bs a}\ \otimes \ \tilde L_{\bs b} \ =\ \sum _{\bs c}\
\tilde N_{\bs a, \bs b, \bs c}\ \tilde L_{\bs c} \  .
\eean
If $\bs a, \bs b, \bs c$ are such that $a_1, b_1, c_1$ are not greater than $d - N$, then
\be
\tilde N_{\bs a, \bs b, \bs c}\ = \ N_{\bs a, \bs b, \bs c} \ ,
\ee
see \cite{M}.

Any $gl_{N+1}$-module $\tilde L$ has a natural structure of an $sl_{N+1}$-module denoted $L$.

The  Littlewood-Richardson rule and the  above remarks imply the following corollary.

\begin{cor}\label{multiplicity}
The number of $V\in Gr (N+1, \C_d[x])$ with prescribed ramification conditions at 
$z_1, \dots , z_n$, such that \Ref{ram cond} is satisfied, is not greater than the multiplicity 
of the trivial $sl_{N+1}$-module in the tensor product of $sl_{N+1}$-modules
\be
 L_{\bs a (z_1)} \otimes \dots  \otimes
L_{\bs a (z_n)} \otimes  L_{\bs a (\infty)} \ .
\ee
\end{cor}

Conjecturally, for almost all $z_1, \dots z_n$ the number of such $V$
is equal to the above multiplicity. For $N=1$ the statement of this conjecture is one of corollaries 
of the main result in \cite{ScV}.

\section{The case of the Lie algebra $sl_{N+1}$}

In this section we study populations of critical points associated with
$sl_{N+1}$. Let  $\al_i, \ i = 1, \dots , N$, be simple roots of $sl_{N+1}$.
 We have
$( \al_i, \al_i ) = 2$ and $(\al_i, \al_{i\pm 1}) = -1$ with all other scalar 
products equal to zero.

\subsection{Reproduction procedure for $\g = sl_2$ 
and associated differential operators.}

Let $\La_s, s = 1, \dots , n,$ be dominant integral $sl_2$-weights.
Let $z_1, \dots , z_n\ \in \ \C$, \
$l\in\Z_{\geq 0}$, and $\Lambda_\infty = \sum_s \Lambda_s - l\alpha$.

Let a polynomial $y$ represent a critical point of the master function 
$\Phi(\bs t; \bs z; \bs \Lambda, \Lambda_\infty)$.

For $\g=sl_2$,  system (\ref{wronskian-critical eqn})
is a single equation $ W(y, \tilde{y}) = T$, where $T = \prod_{s}(x-z_s)^{(\La_s, \al)}$.

Let $\tilde y$ be a solution of this equation.
Let $V$ be the complex two-dimensional vector space of polynomials 
spanned by $y$ and $\tilde y$. 
The vector space $V$ is the space of solutions of the differential equation
\be
F \ u''\ - \ G\ u' +\ H\ u\ =\ 0\ ,
\ee
where 
\be
 F(x) = \prod_{s=1}^n (x - z_s), \ G(x) = F(x) \ln'(T(x)), \ H(x) 
= (G(x) y'(x) - F(x) y''(x))/ y(x).
\ee
\begin{lem}\label{sl2 pop}
The projective line $\bs P(V)$ is the population of critical points associated with 
$sl_2$, dominant integral weights 
$\Lambda_1, \dots , \Lambda_n$, points $z_1, \dots , z_n$, and originated at $y$.
Almost all points of the line represent critical points. The exceptions form a
finite subset of the line.
\hfill $\square $
\end{lem}

For any integer $k \geq 1$ it is useful to  introduce a linear differential operator of order $k+1$
whose kernel is spanned by polynomials $y^i \tilde y^{k-i}$, $i=0,\dots,k$. 
Namely, for any $g\in V$, set
\bea
D_k(g) &=& (\partial + k \ln' g - k \ln' T))
\ \dots \ (\partial
+ (-k+2) \ln' g - \ln' T)\ (\partial - k \ln' g)  \notag\\
&=& \prod^{0\ \to \ k}_{i}(\partial + (k-2i) \ln' g + (i-k) \ln' T)\ .
\eea


\begin{lem}\label{sl2 prop}
The operator $D_k(g)$ does not depend on the choice of a non-zero element $g \in V$, 
and the kernel of  $D_k(g)$ is spanned by polynomials $y^i \tilde y^{k-i}$,\
$i = 0, \dots , k$. 
\end{lem}
\begin{proof} 
Let $g, \tilde g$ be a basis in $V$.
By induction on $i$, we prove   that the kernel of the product of
$i$ right linear factors
in $D_k(g)$ is spanned by functions
$g^k, g^{k-1} \tilde g, \dots , g^{k-i+1} \tilde g^{i-1}$. 

For $i=1$ the statement is clear. Suppose that the statement  is proved for
$i<i_0$. Then the kernel of the product of $i_0$ right factors is an $i_0$-dimensional 
space which contains $g^k,$ $ g^{k-1} \tilde g,$ 
$ \dots,$ $ g^{k-i_0+2} \tilde g^{i_0-2}$. Note that the
expansion of the product of $i_0$ right factors starts with
\be
\partial^{i_0}\ -\  \left( i_0(k-i_0+1) \ln'y\ +\
i_0(i_0-1)/2 \ln' T \right) \partial^{i_0-1}\ +\ \dots \;.  
\ee
Therefore it is enough to show that 
\be
W(g^k, g^{k-1} \tilde g, \dots , g^{k-i_0+1} \tilde g^{i_0-1} )\ 
=\ c\  g^{i_0(k-i_0+1)}\ T^{i_0(i_0-1)/2} 
\ee
for some non-zero constant $c$. This readily follows from Lemmas
\ref{f wronskian} and \ref{fg id}.
\end{proof}

\subsection{The differential operator 
associated to  an $sl_{N+1}$-population}

Let $\La_s, s = 1,  \dots , n,$ be dominant integral $sl_{N+1}$-weights.
Let $z_1, \dots , z_n \ \in \ \C$.
Let $l_i\in\Z_{\geq 0},\ i = 1, \dots  , N$, and 
$\Lambda_\infty = \sum_{s=1}^n \Lambda_s - \sum_{i=1}^N l_i\alpha_i$.
Let the polynomials $T_i,\ i = 1, \dots , N$, be defined by \Ref{T}.

Let $\bs y =(y_1, \dots, y_N)$ be an $N$-tuple of non-zero polynomials.
We adopt the notation $y_0 = y_{N+1} = 1$. For any $\bs y$,
define a  linear differential operator of order $N+1$ 
with meromorphic coefficients 
\bean\label{Diff}
D ( \bs y) & =&  (   \partial  - 
\ln'  ( \frac { \prod_{s=1}^N T_s } { y_N }  ) )\
( \partial - \ln' ( \frac {y_N  \prod_{s=1}^{N-1} T_s} {y_{N-1} }  ) )
\dots ( \partial  -  \ln'  ( \frac { y_2 T_1 }{ y_1 }  ) ) \
( \partial - \ln' (  y_1 ) ) 
\notag
\\
& = &
 \prod^{0\ \to\ N}_i\ (\
\partial \ -\ \ln' \ ( \
\frac{ y_{N+1-i}  \prod_{s=1}^{N-i} T_s }{ y_{N-i} }\ ) \ )\ .
\eean
Notice that the first coordinate $y_1$ of the $N$-tuple  belongs to 
the kernel of the operator $D(\bs y)$.

\begin{theorem}\label{ind thm}
Let $P$ be an $sl_{N+1}$-population of critical points associated with weights
$\Lambda_1, \dots  , \Lambda_n$, points $z_1, \dots  , z_n$, and originated at some $N$-tuple
$\bs y^0$. Then the differential operator $D(\bs y)$ does not depend on the choice
of $\bs y$ in $P$.
\end{theorem}

This differential operator is called {\it associated with the population} and is denoted $D_P$.

\begin{proof}
First we prove that if $\bs y = (y_1, \dots  , y_N)$ and 
$\tilde{\bs y} = (\tilde y_1, \dots , \tilde y_N)$ are members of the population 
which represent critical points, then $D ( \bs y )\ =\ D ( \tilde {\bs y} )$.
To show that statement,
 it is sufficient to consider the case, when $\bs y$ and 
$\tilde{\bs y}$ differ only  in one coordinate. Namely, we assume that there is an index $i$ such
that $y_j\ =\ \tilde y_j$, for all $j, \ j \neq  i$,  and
$W(y_i, \tilde y_i) \ = \ T_i\  y_{i-1}\ y_{i+1}$. 
In this case, all linear factors of $D(\bs y)$ and
$D(\tilde{\bs y})$ are the same except the two linear factors 
which involve $y_i$ or $\tilde y_i$. So to prove 
that $D(\bs y)\ =\ D(\tilde{\bs y})$ we have to show that for any function $u(x)$ we have
\bea
( \partial\  - \ \ln' ( \frac{y_{i+1} \ \prod_{s=1}^{i}{T_s}}{y_{i}}))
 \ ( \partial\  -\  \ln' ( \frac{y_i \ \prod_{s=1}^{i-1}{T_s}}{y_{i-1}} ) )\ u\ =
\\
( \partial\ -\ \ln' ( \frac{ y_{i+1}\  \prod_{s=1}^{i}{T_s}}{\tilde y_{i}} ) )\
( \partial\ -\ \ln' ( \frac{\tilde y_i\  \prod_{s=1}^{i-1}{T_s}}{y_{i-1}} ) )\ u .
\eea
After the change of variables $v\ =\ u\ (\prod_{s=1}^{i-1}T_s)\ / \ y_{i-1}$, we have to prove that
 \bea
( \partial\ -\ \ln'( \frac{y_{i+1}y_{i-1}T_i}{y_{i}} ) ) \ ( \partial\ -\ \ln' ( y_i ) )\ v\
=\ ( \partial\ -\ \ln' ( \frac{ y_{i+1}y_{i-1}T_i}{\tilde y_{i}} ) )\ 
( \partial\ -\ \ln' ( \tilde y_i ) )\ v\ .
\eea
This identity  follows from the $sl_2$ case, see 
Lemma \ref{sl2 prop} with $k=1$.

The $N$-tuples $\bs y$, which represent critical points, form a dense subset in the population.
The operator $D ( \bs y )$ continuously depends on $\bs y$ and is constant on a dense subset 
of the population. Hence it is constant on the whole  population.
\end{proof}

\begin{cor}\label{first coor cor}
Let $\bs y$ be a member of an $sl_{N+1}$-population $P$. 
Then the first coordinate $y_1$ of ${\bs y}$ lies in the kernel of the operator $D_P$.
\end{cor}

\subsection{Equation $D_P u = 0$ is Fuchsian} \label{all polynom section}
In this section we show that all solutions of equation 
\bean\label{DP eq}
D_P \ u\ =\ 0
\eean
 are polynomials,  and hence the equation is Fuchsian.

Let $\bs y$ be a member of the population $P$. Assume that $\bs y$ represents a critical point. 
System (\ref{wronskian-critical eqn}) for $\g = sl_{N+1}$  takes the form
\bean\label{wronskian-critical N}
W(y_i,\tilde{y}_i)\ = \ T_i \ y_{i-1} \ y_{i+1}\ , \qquad i = 1, \dots , N\ .
\eean

We start with a construction of polynomials $u_1, \dots , u_{N+1}$,
satisfying equation \Ref{DP eq}. 

Set $u_1=y_1$. 
The polynomial $u_1$ is a solution of \Ref{DP eq}  by Theorem \ref{ind thm}.

Let $u_2$ be a polynomial such that
$W(u_1,u_2) \ = \ T_1 y_2$. Such a polynomial exists by Lemma
\ref{wronskian-critical lem}.
The polynomial $u_2$ is a solution of \Ref{DP eq}  by Theorem \ref{ind thm}.

Let $\tilde y_2$ be a polynomial such that $W ( y_2, \tilde y_2 ) \ = \ T_2 y_1 y_3$. 
Such a polynomial exists by Lemma
\ref{wronskian-critical lem} and  can be chosen so that
$( y_1, \tilde y_2 , \dots , y_N )$ is generic and therefore
represents a critical point, see Lemma \ref{sl2 pop}. 
Choose a polynomial $u_3$ to satisfy  equation
$W(u_3,y_1) \ =\ T_1 \tilde y_2 $. 
The polynomial $u_3$ is a solution of \Ref{DP eq}  by Theorem \ref{ind thm}.

In general, to construct a polynomial
$u_{i+1}$ we find $\tilde {\bs
  y}_i=(y_1,\dots,\tilde y_{i},\dots,y_N)$, such that $\tilde {\bs y}_i$ is
  generic  and $W(y_i,\tilde y_i)=T_iy_{i-1}y_{i+1}$ and then repeat
  the construction for $u_i$ using $\tilde{\bs y}_i$ instead of $\bs y$.

Let $V$ be the complex vector space spanned by $u_1,\dots,
u_{N+1}$.  We show below that the space $V$ is the kernel of the operator $D_P$.
Thus the space $V$ does not depend on  the choice of $\bs y \in P$.
We call the space {\it the fundamental space of the population} $P$
and denote $V_P$.

\begin{lem}\label{wr(u) lem}
For $i=1,\dots,N+1$, we have 
\be
W(u_1, \dots, u_i)\ =\ y_i \ T_1^{i-1}\  T_2^{i-2}\  \dots\  T_{i-1}\ .
\ee
\end{lem}
\begin{proof}
We proceed by induction on $i$. 
The case $i=1$ is trivial. Suppose the statement
is proved for $i=1, \dots , i_0-1$. By induction hypothesis,
\bea
W(W(u_1, \dots , u_{i_0-1}), W(u_1, \dots , u_{i_0-2}, u_{i_0} ) )\ &= &\
W(y_{i_0-1}\prod_{s=1}^{i_0-2}T_s^{i_0-1-s}, \ \tilde
y_{i_0-1} \prod_{s=1}^{i_0-2}T_s^{i_0-1-s} )\ \\
=\ (\prod_{s=1}^{i_0-2}T_s^{i_0-1-s} )^2\ W(y_{i_0-1},\tilde y_{i_0-1})\ & =&
\ (\prod_{s=1}^{i_0-2} T_s^{i_0-1-s} )^2 \ y_{i_0-2}\  y_{i_0}\ T_{i_0-1}\ .
\eea
On the other hand, by Lemma \ref{wr id 1} we have
\be
W(W(u_1,\dots,u_{i_0-1}), \ W(u_1,\dots,u_{i_0-2},u_{i_0}))\ = \
W(u_1,\dots,u_{i_0})\ W(u_1,\dots,u_{i_0-2})\ .
\ee
That implies the lemma.
\end{proof}

\begin{lem}\label{kernel prop} The space $V$ is the kernel of the operator $D_P$.
\end{lem}
\begin{proof}
The polynomials $u_1,\dots, u_{N+1}$ are linearly independent solutions
of equation \Ref{DP eq}.
\end{proof}

\subsection{Singular points and exponents of the equation $D_P u = 0$}\label{no b}
Choose a basis $u_1, \dots , u_{N+1}$ in $V_P$ as in Section \ref{all polynom section}.

\begin{lem}
Singular points of $D_P u = 0$ lie in 
$\{ z_1, \dots , z_n, \infty \}$.
\end{lem}
\begin{proof}
The operator $W (u_1, \dots , u_{N+1})\, D_P$ has polynomial coefficients.
By Lemma \ref{wr(u) lem} zeros of $W ( u_1, \dots , u_{N+1} )$ 
lie in $\{ z_1, \dots , z_n\}$.
\end{proof}

\begin{lem}\label{z exp prop}
For $s = 1, \dots , n$, the  exponents of  $D_P u = 0$ at $z_s$ are
\be
 \ 0,\ (\La_s + \rho, \ \al_1) ,\ (\La_s + \rho, \ \al_1 + \al_2) , \dots ,
(\La_s + \rho, \ \al_1 + \dots + \al_N)\ .
\ee
\end{lem}
\begin{proof}
For $i = 1, \dots , N+1$, let $V_i$ be the space spanned by
$u_1, \dots , u_i$. Let $D_i$ be
a linear differential operator of order $i$ whose kernel is $V_i$.
Such an operator is unique up to multiplication by a function.
In particular, we have $V_{N+1} = V_P$ and $D_{N+1} = D_P$.

We prove by induction on $i$, that the exponents of $D_i$ at $z_s$
are
\be
\ 0,\ (\La_s + \rho, \ \al_1) ,\ (\La_s + \rho, \ \al_1 + \al_2) , \dots ,
(\La_s + \rho, \ \al_1 + \dots + \al_i)  \ .
\ee
For $i=1$ we have $u_1 (z_s) \neq 0$, hence the exponent of $D_1$ at $z_s$
is zero. 

Suppose we proved the statement for $i = 1, \dots , i_0 - 1$.
By Lemma \ref{wr(u) lem} the Wronskian $W(u_1, \dots , u_{i_0})$ has zero at
$z_s$ of order $( \La_s,\ \sum_{j=1}^{i_0-1} (i_0-j) \al_j)$. That implies our
statement for $i = i_0$.
\end{proof}

\begin{cor}\label{no base}
The vector space $V_P$ has no base points.
\end{cor}

By Theorem \ref{inf weight thm} there exists $\tilde{\bs y} \in P$ such
that the weight of $\tilde{\bs y}$ at infinity with respect to 
$\La_1, \dots , \La_n$  is integral dominant.
Remind that the weight at infinity is 
$\tilde \La_\infty \ =\ \sum_{s=1}^n \La_s - \sum_{i=1}^N \tilde l_i\al_i$, with
$\tilde l_i = \deg \tilde y_i$.

\begin{lem}\label{inf exp prop}
The exponents of $D_P u = 0$ at infinity are
\be
 \tilde l_1,\ \tilde l_1 + (\tilde \La_\infty + \rho,\ \al_1),\
\tilde l_1 + (\tilde\La_\infty + \rho,\ \al_1+\al_2), \ \dots ,\
\tilde l_1 + (\tilde\La_\infty + \rho,\ \al_1 + \dots + \al_N) \ .
\ee
\end{lem}
The proof is similar to the proof of Lemma \ref{z exp prop}.

\begin{cor}\label{ram points}
The ramification points of $V_P$ are $z_1, \dots , z_n, \infty$. 
For $s = 1, \dots , n$, the ramification condition $\bs a(z_s)$
has the form $a_i (z_s) = \sum_{j=1}^{N+1-j} (\La_s, \al_j)$,\
$i = 1, \dots , N+1$.
The ramification condition $\bs a(\infty)$ has the form
\be
a_i (\infty)  =  d - N - \tilde l_1 - (\tilde \La_\infty, \ \al_1 + \dots + \al_{i-1})\ ,
\ee
$i = 1, \dots , N+1$.
\end{cor}

\subsection{Vector spaces of polynomials and critical points}\label{ha-ha}
Let $V \in Gr(N+1, \C_d[x])$ be an $N+1$-dimensional vector subspace
in $\C_d[x]$ with no base points. 
Assume that the ramification set of space $V$ consists of  points 
$z_1, \dots , z_n, \infty$
with ramification conditions $\bs a(z_1), \dots , \bs a(z_n), \bs a(\infty)$, respectively.

Introduce dominant integral $sl_{N+1}$-weights $\La_1, \dots , \La_n, \tilde \La_\infty$.
These weights will be called {\it the weights of $V$ }.

For $s = 1, \dots , n$, let $\La_s$ be defined by  conditions
$(\La_s, \al_i) = a_{N+1-i}(z_s) - a_{N+2-i}(z_s)$, \ $i = 1 , \dots , N$.
Let $\tilde \La_\infty$ be defined by  conditions
$(\tilde \La_\infty, \al_i) = a_{i}(\infty) - a_{i+1}(\infty)$, \ $i = 1 , \dots , N$.

Let $FL(V)$ be the variety of all full flags
\be
\mathcal F \ = \ \{ 0 \subset F_1 \subset F_2 \subset \dots \subset
F_{N+1} = V \}\ 
\ee 
in $V$. For any $\mathcal F \in FL(V)$ define an $N$-tuple of polynomials
$\bs y^{\mathcal F} =  (y_1^{\mathcal F}, \dots , y_N^{\mathcal F})$\ 
as follows. Let $u_1, \dots , u_{N+1}$ be a basis in $V$ such that for any $i$ the polynomials
$u_1, \dots , u_i$ form a basis in $F_i$. 
Say that this {\it basis is adjusted } to the flag $\mc F$ and the flag
$\mc F$ {\it is generated} by the basis $u_1, \dots , u_{N+1}$. We set
\bean\label{div}
y_i^{\mathcal F} \ = \  W(u_1, \dots , u_i) \ T_1^{1-i}\ T_2^{2-i} \ \dots\ T_{i-1}^{-1}\ ,
\eean
where polynomials $T_1, \dots , T_N$ are given by \Ref{T}. 
The polynomials $y_1^{\mathcal F}, \dots , y_N^{\mathcal F}$  are uniquely
determined by the flag $\mathcal F$ up to multiplication of polynomials by non-zero numbers.

This gives a morphism
\bean\label{morphism}
\beta \ : \ FL(V) \ \to \PCN , \qquad
\mathcal F \ \mapsto \ y^{\mathcal F}\ ,
\eean
which will be called {\it the generating morphism of $V$}.

\begin{theorem}\label{pol-crit}
${ \ }$

\begin{enumerate}
\item[1.] The generating morphism  defines an isomorphism of \ $FL(V)$
and its image \newline $\beta (FL(V)) \subset \PCN$.
\item[2.]
Let $O$ be the orbit of $\tilde \La_\infty$ with respect to the 
shifted action of the Weyl group. Then
there exist a weight $\La_\infty \in O$ and a critical point $\bs y^0$
of the master function $\Phi (\bs t; \bs z; \bs \La, \La_\infty)$, such that
$\beta (FL(V))$ is the population of critical points originated at $\bs y^0$.

\item[3.]
The initial vector space $V$ is the fundamental space of the population $\beta (FL(V))$.

\item[4.]
Every $sl_{N+1}$-population of critical points is given by this construction
starting from a suitable space $V$.
\end{enumerate}
\end{theorem}

Theorem \ref{pol-crit} is proved in Section \ref{Proof}.

Theorem \ref{pol-crit} and corollaries of Section \ref{no b} say that the 
$sl_{N+1}$-populations of critical points associated with prescribed highest weights at
points $z_1, \dots z_n, \infty$ are in one-to-one correspondence with
elements $V \in Gr(N+1, \C_d[x])$ with no base points and prescribed ramifications at
$z_1, \dots , z_n, \infty$.


\subsection{The upper bound for the number of populations}
Let $V$ be a space of polynomials of dimension $N+1$ with no base points.
Let $z_1,\dots, z_n, \infty$ be its ramification points. There are three
equivalent ways to discribe ramification properties of $V$:
\begin{enumerate}
\item[$\bullet$] indicating ramification conditions $\bs a(z_1),\dots,\bs a(z_n),\bs
a(\infty)$ of $V$ in $\C_d[x]$, where $d$ is large enough to have
$V\subset\C_d[x]$,
\item[$\bullet$] indicating  exponents $\bs m(z_1),\dots,\bs m(z_n), \bs m(\infty)$ of the
linear differentail equation $Du=0$ with kernel $V$,
\item[$\bullet$] indicating integral dominant $sl_{N+1}$-weights
$\La_1,\dots,\La_n,\La_\infty$.
\end{enumerate}
According to Sections \ref{no b} and \ref{ha-ha} the relations between these data are
\bea
(\La_s,\al_i) = a_{N+1-i}(z_s) - a_{N+2-i}(z_s), 
\qquad && 
(\La_\infty,\al_i) = a_{i}(\infty) - a_{i+1}(\infty),
\\
a_i(z_s) = (\La_s, \al_1 + \dots + \al_{N+1-j}),
\qquad && 
a_i(\infty) = d - N - l_1 - (\La_\infty, \al_1 + \dots + \al_{i-1}), 
\\
m_i(z_s) = a_{N+2-i} + i - 1,
\qquad && 
m_i(\infty) = d - N - a_i(\infty) + i - 1.
\eea
Here $l_1 = (\sum_{s=1}^{n}\La_s - \La_\infty, \om_1)$ where $\om_1$ is the
first fundamental weight, $(\al_i,\om_1) = \delta_{1i}$ for all $i$.

Denote $L_\La$ the irreducible $sl_{N+1}$-module with highest weight $\La$.
Let $w$ be an element of the Weyl group.
Combining Corollary \ref{ram points} and Corollary \ref{multiplicity} we get the following statement.

\begin{theorem}\label{estimate}
The number of $sl_{N+1}$-populations of 
critical points which contain a critical 
point of the master function $\Phi(\bs t;\bs z;\bs \La, \ w\cdot
\La_\infty)$  is not greater than
the multiplicity of the  $sl_{N+1}$-module $L_{ \La_\infty}$
in the tensor product of  $sl_{N+1}$-modules
$ L_{\La_1} \otimes \dots  \otimes L_{\La_n}$.
\end{theorem}

We expect that for almost all $z_1, \dots, z_n$ the number of populations indicated in Theorem
\ref{estimate} is equal to the above multiplicity.

\subsection{Proof of Theorem \ref{pol-crit} }\label{Proof}

\begin{lem}\label{diff}
The differential of $\beta$ is non-degenerate at any point 
$\mathcal F \in FL(V)$.
\end{lem}

\begin{proof}
Let $u_1, \dots , u_{N+1}$ be a basis in $V$ adjusted to $\mc F$.
For $k = 1, \dots , N+1$, and a subset $\bs i=(i_1,\dots,i_k) \subset \{1, \dots , N+1\}$,
 define a polynomial 
\be
u_{k,\bs i}\ =\ W(u_{i_1}, \dots , u_{i_k} )\ T_1^{1-i} \ T_2^{2-i}\dots T_{i-1}^{-1}\ .
\ee
For a fixed $k$, the polynomials $u_{k,(1, 2, \dots , k-1 , j)}$, $j = k+1, \dots , N+1$, are
linearly independent.

In the neighborhood of $\mathcal F$ in $FL(V)$, every flag is
generated by a basis $\tilde u_1, \dots , \tilde u_{N+1}$ of the form
$\tilde u_i=  u_i + \sum_{j,\ j>i} c_{ij} u_j$. The coefficients 
$c_{ij}$ are local coordinates on $FL(V)$ at $\mc F$. 

Consider a curve $\mc F^s$ through $\mc F$, where the flag $\mc F^s$ is generated by a basis
\be
\tilde u_i^s\ =\ u_i\ +\ s\ \sum_{j,\ j>i}\ c_{ij}\ u_j, \qquad i = 1, \dots , N+1\ .
\ee
The tangent vector to this curve at $s = 0$ has coordinates $\{ c_{ij}\}_{j>i}$. 
The image of the tangent vector under the differential of $\beta$ is
the non-zero vector
\be
(\ \sum_{j>1}\ c_{1j} \ u_j \ , \  \dots \ ,\ \sum_{j>k}\ c_{kj}\ u_{k,\ (1, \dots, k-1, j)}\ ,
\  \dots\ , \ c_{N,N+1}\ u_{N, \ ( 1, 2, \dots N-1, N+1)}\ ) \ .
\ee
\end{proof}

Let $\mc F \in FL(V)$ and let $u_1, \dots , u_{N+1}$ be
a basis in $V$ adjusted to $\mc F$.
For $i \in \{1, \dots , N\}$ and a non-zero number $c$, consider the flag 
\be
\mc F^{(i)} \ = \ \{ 0 \subset F_1 \subset \dots F_{i-1} \subset \tilde F_i
\subset F_{i+1} \subset \dots \subset F_{N+1}  \}\ 
\ee
generated by the basis $u_1, \dots , u_{i-1}, u_i + c u_{i+1}, u_{i+1}, \dots , u_{N+1}$.

\begin{lem}\label{Imm}
The image of $\beta$ consists of fertile $N$-tuples with respect to polynomials
$T_1, \dots , T_N$. The $N$-tuple $\beta (\mc F^{(i)})$
is an immediate descendant of $\beta(\mc F)$ in the direction $i$.
\end{lem}
\begin{proof} We have $y^{\mc F^{(i)}}_j  =  y^{\mc F}_j$ 
for $j\ \neq i$ and
\bea
W( y^{\mc F}_i ,\ y^{ {\mc F}^{(i)} }_i ) & = &
c\ W ( W (u_1 , \dots , u_{i}), \ W ( u_1, \dots , u_{i-1}, u_{i+1}))\
(\prod_{j=1}^{i-1}T_j^{i-j} )^2 
\\
&=& c \ W ( u_1, \dots , u_{i-1} )\ W(u_1, \dots , u_{i+1} )\
( \prod_{j=1}^{i-1} T_j^{i-s} )^2
\\
&=& c \ T_i\ y^{\mc F}_{i-1}\ y^{\mc F}_{i+1}\ .
\eea
 \end{proof}

\begin{lem} The generating morphism $\beta$ maps distinct points to distinct.
\end{lem}
\begin{proof} If $\bs y = \beta (\mc F) = \beta (\tilde {\mc F})$, then $\mc F = \tilde{\mc F}$.
Indeed, let  $\mc F$ be generated by $u_1, \dots , u_{N+1}$ and $\tilde{\mc F}$ by
$\tilde u_1, \dots , \tilde u_{N+1}$. Clearly, $u_1$ and $\tilde u_1$ are proportional. 
Since $ W ( u_1, u_2 ) = c\, T_1 y_2$
and $ W ( \tilde u_1, \tilde u_2 ) = \tilde c\, T_1 y_2$ for suitable non-zero $c$ and $ \tilde c$,
the span of $u_1, u_2$ is equal to the span of $\tilde u_1, \tilde u_2$, and so on.
\end{proof}

Part 1 of Theorem \ref{pol-crit} is proved.


\begin{prop} \label{generic prop}
There exists $\mc F \in FL(V)$ such that $\beta (\mc F)$ is generic with respect to
integral dominant weights $\La_1, \dots, \La_n$ of 
the Lie algebra $sl_{N+1}$ and points $z_1, \dots , z_n$.
\end{prop}
\begin{proof}
Let $\bs y = (y_1, \dots , y_N) \in \beta(FL(V))$. The image of $\beta$ consists of fertile points.
Hence if $y_i$ has a multiple root, then the multiple root is either a root of $y_{i-1}$ or
a root of $y_{i+1}$ or lies in  $\{ z_1, \dots , z_n\}$.

\begin{lem}
For almost all $\bs y \in \beta (FL(V))$ we have $y_i(z_s) \neq 0$ for all $i$ and $ s$.
\end{lem}
\begin{proof}
Let $\mc F \in FL(V)$, $
\mathcal F \, = \, \{ 0 \subset F_1 \subset F_2 \subset \dots \subset
F_{N+1} = V \}$.
For $i \in \{1, \dots , N\}$, let  $D_i$ be a linear differential operator whose kernel 
is $F_i$. The operator $D_i$ is defined up to multiplication by a function and
is Fuchsian. Let $s \in \{1, \dots , n\}$. 
For almost all $\mc F$ and every $i$, 
the exponents of $D_i$ at $z_s$ are 
\be
0,\ (\La_s + \rho, \ \al_1) ,\ (\La_s + \rho, \ \al_1 + \al_2) , \dots ,
(\La_s + \rho, \ \al_1 + \dots + \al_{i-1})\ .
\ee
For such $\mc F$, clearly $y_1^{\mc F}, \dots , y_N^{\mc F}$ do not have roots at
$z_s$.
\end{proof}

\begin{lem}
For every $i$ and almost all $\bs y \in \beta (FL(V))$, the polynomials
$y_i$ and $y_{i+1}$ do not have common roots.
\end{lem}
\begin{proof}
The set $X_i$ of all $\mc F \in FL(V)$, such that $y_i^{\mc F}$ and $y_{i+1}^{\mc F}$
have common roots, is an algebraic subset of $FL(V)$.
For every $z \in \C$, the set of all $\mc F \in FL(V)$, such that $y_i^{\mc F}(z) = y_{i+1}^{\mc F}(z)
= 0$, has codimension two in $FL(V)$ . Hence $X_i$ is a proper algebraic subset of $FL(V)$.
\end{proof}
The above lemmas prove the proposition.
\end{proof}

Let $\bs y^0 \in \beta(FL(V))$ be a generic point. By Corollary  \ref{fertile cor},
$\bs y^0$ represents a critical point of an $sl_{N+1}$-master function, associated with 
integral dominant weights $\La_1, \dots , \La_n$ ( 
which are defined at the beginning of this section ), 
suitable integral weight at infinity $\La_\infty$, complex numbers $z_1, \dots , z_n$.

\begin{lem}
The set $\beta(FL(V))$ is the population of critical points originated at $\bs y^0$.
Moreover, the fundamental space of the population  $\beta(FL(V))$ is the initial space $V$.
\end{lem}
\begin{proof}
Let $P$ be the population of critical points originated at $\bs y^0$. Lemma \ref{Imm}
shows that $P$ is contained in the image of $\beta$. It is easy to see that
$P = \beta (FL(V))$. The second statement of the lemma directly
follows from the definition of the generating morphism $\beta$.
\end{proof}
To prove part 2 of Theorem \ref{pol-crit} it is enough to prove the following lemma.
\begin{lem}\label{iNf}
The weight $\La_\infty$ lies in the orbit $O$ of the weight $\tilde \La_\infty$.
\end{lem}
\begin{proof}
By the definition in Section \ref{FER}, the weight at infinity of an $N$-tuple $\bar {\bs y}$
with respect to weights $\La_1, \dots , \La_n$ is 
$\bar \La_\infty \ =\ \sum_{s=1}^n \La_s - \sum_{i=1}^N \bar l_i\al_i$, with
$\bar l_i = \deg \bar y_i$.

To prove the lemma it is enough to find  $\bar {\bs y} \in \beta(FL(V))$ whose weight at infinity 
is $\tilde \La_\infty$. 

The space $V$ has $\bs a(\infty)$ as the ramification index at infinity.
Hence $V$ has a basis of polynomials $u_1, \dots , u_{N+1}$ 
such that deg $u_i = d - N - a_i(\infty) + i - 1$.
 Let $\bar {\mc F} \in FL(V)$ be the flag generated by this basis and
let  $\bar {\bs y} = \beta (\bar{\mc F})$. 

It is clear that  for any $\bs y \in \beta(FL(V))$, we have $\deg\,y_i \geq  \deg\,\bar {y}_i$,
$i = 1, \dots , N$. Therefore by Theorem \ref{inf weight thm}, the weight 
$\bar \La_\infty \ =\ \sum_{s=1}^n \La_s - \sum_{i=1}^N \bar l_i\al_i$, with
$\bar l_i = \deg \bar y_i$, is the dominant weight lying in the orbit $O$.

We show  that  $\tilde \La_\infty = \bar \La_\infty$.
Indeed, set $\tilde l_1 = d - N - a_1(\infty)$. Then,  according to our definition of $\tilde \La_\infty$,
the exponents at infinity of  $V$ are given by the formula of Lemma \ref{inf exp prop}. 

At the same time, by Lemma \ref{inf exp prop}, the exponents at infinity of $V$ are
\be
 \bar  l_1,\ \bar l_1 + (\bar \La_\infty + \rho,\ \al_1),\
\bar l_1 + (\bar\La_\infty + \rho,\ \al_1+\al_2), \ \dots ,\
\bar l_1 + (\bar\La_\infty + \rho,\ \al_1 + \dots + \al_N) \ .
\ee
Comparing the two formulas, we get  $\tilde \La_\infty = \bar \La_\infty$.
\end{proof}
The previous lemmas prove parts 2 and 3 of Theorem \ref{pol-crit}. Part 4 of the theorem follows 
from Lemma \ref {wr(u) lem}. Theorem \ref{pol-crit} is proved.

\subsection{Populations and Bruhat cells}\label{Bruhat section} 
Let $V$ be a complex vector space of dimension $N+1$.
Fix a full flag $\mc F^0 \in FL(V)$. Then there is 
a Schubert-type cell decomposition of the variety $FL(V)$ of
full flags in $V$. 

 For any $\mc F \in FL(V)$ define a permutation
$w(\mc F)$ in the symmetric group $S^{N+1}$ as follows. 
Define $w_1(\mc F)$ as the minimum of $i$ 
such that  $F_1 \subset F^0_i$. Fix a basis vector $u_1 \in F_1$. 
Define $w_2(\mc F)$ as the minimum of $i$ such that  
there is a basis in $F_2$ of the form $u_1, u_2$ with $u_2 \in F_i^0$.
Assume that 
$w_1(\mc F), \dots , w_j(\mc F)$ and $u_1,\dots,u_j$ are determined.
Define $w_{j+1}(\mc F)$ as the minimum of $i$ such that  
there is a basis in $F_{j+1}$ of the form $u_1, \dots , u_j, u_{j+1}$ with 
$u_{j+1} \in F_i^0$. As a result of this 
procedure we define $w(\mc F) = (w_1(\mc F), \dots , w_{N+1}(\mc F))\ \in S^{N+1}$
and a basis $u_1, \dots , u_{N+1}$ which generates $\mc F$ and such that
$u_i \in F^0_{w_i(\mc F)}$.

For $w \in S^{N+1}$, define 
\be
G^{\mc F^0}_w \ =\ \{\ \mc F \in FL(V)\ , \ w(\mc F) \,=\, w\ \}\ .
\ee 
The algebraic variety $G^{\mc F^0}_w$ is
called {\it the Bruhat cell } associated with $\mc F^0$ and  $w \in S^{N+1}$ .
The set of all Bruhat cells form a cell decomposition of $FL(V)$:
\be
FL(V)=\sqcup_{w\in S^{N+1}}G^{\mc F^0}_w.
\ee

Let $V \in Gr ( N+1, \C_d[x] )$ be an $N+1$-dimensional space of polynomials 
with ramification points  $z_1, \dots , z_n, \infty$ and ramification
conditions 
$\bs a(z_1), \dots , \bs a(z_n), \bs a(\infty)$. Let $\La_1, \dots
\La_n, \tilde \La_\infty$ be the corresponding integral dominant
$sl_{N+1}$-weights. 

The space $V$ has a distinguished flag $\mc F^\infty \in FL(V)$, where
$F^\infty_i$ consists of all polynomials in $V$ of degree not greater than 
$d - N - a_i(\infty) + i - 1$.
Then the variety $FL(V)$ has the corresponding
Bruhat cell decomposition $FL(V)=\sqcup_{w\in S^{N+1}}G_w^{\mc F^\infty}$.

For any $N$-tuple $\bs y$ of polynomials, introduce a vector $\bs l(\bs y) = 
(l_1(\bs y), \dots , l_N(\bs y))$, where $l_i(\bs y) = \deg y_i$.

Recall that the generating morphism $\beta : FL(V) \to \PCN $ is given by \Ref{morphism}.
By definition, any two flags $F_1$ and $F_2$ belong to the same Bruhat cell $G_w^{\mc
  F^{\infty}}$ if and only if
$\bs l(\beta(F_1))=\bs l(\beta(F_2))$. Set $\bs l^w=\bs l(\beta(F))$
for $F\in G_w^{\mc F^\infty}$.

We identify the Weyl group $W$ of $sl_{N+1}$ with the symmetric group $S^{N+1}$ 
sending the generating reflections $s_i$ to the simple transpositions $(i, i+1)$.

\begin{lem}\label{sln bruhat lem} For $w \in S^{N+1}$ we have 
\be
\sum_{i=1}^{N} \ l_i^w\ \al_i\ =\ \sum_{s=1}^n {\La_s} \ -\ w\cdot \La_\infty\ .
\ee
\end{lem}
\begin{proof} 

Denote $d_i = d - N - a_i(\infty) + i - 1$.
For any permutation $w \in S^{N+1}$ and any $i$, we have the equation 
\be
(w\cdot \La_\infty\ ,\ \al_i)\ =\
d_{w_{i+1}}\ -\ d_{w_i}\ -\ 1 ,
\ee
which can be checked by induction on the length of the permutation $w$.

Also if $w = w(\mc F)$ is the permutation of a flag $\mc F \in FL(V)$,
then $\mc F$ has an adjusted basis $u_1, \dots , u_{N+1}$
such that $\deg u_i = d_{w_i}$ for all $i$. 
Therefore the vector $\bs l^w$ has coordinates
\be
l_i \ =\  \sum _{j=1}^i d_{w_j} \ - \ \sum_{s=1}^n(\La_s+\rho,\ (i-1)\al_1 + \dots + \al_{i-1})\ .
\ee
Combining the two equations we get for any $i$ an equation
\be
2l_i^{w} - l_{i-1}^{w} - l_{i+1}^{w}\ =\ d_{w_i} - d_{w_{i+1}} +
(\sum_{s=1}^n\La_s + \rho, \al_i)\ =\ (\sum \La_s - w\cdot \La_\infty, \al_i) ,
\ee
which is equivalent to the identity of the lemma.
\end{proof}

Let $\La_1, \dots , \La_n,\ \tilde \La_\infty$ be integral dominant $sl_{N+1}$-weights,
$z_1, \dots , z_n$ complex numbers. 
Let $w \in S^{N+1}$. Let $\bs y$ represent a critical point of
the master function  $\Phi(\bs t; \bs z; \bs   \La, w\cdot \tilde \La_\infty)$ and $P$
the population originated at $\bs y$. Let $\beta : FL(V_P) \to \beta (FL(_P))$
be the generating isomorphism.

Let $\C^0_d[x]$ be the affine space of monic polynomials of degree
$d$. Then $\bs y$ determines a point in
$\C^0_{l_1^w}[x]\otimes\dots\otimes\C^0_{l^w_N}[x]$.

\begin{cor} 
The closure in
  $\C^0_{l_1^w}[x]\otimes\dots\otimes\C^0_{l_N^w}[x]$  
of $N$-tuples of polynomials in the population $P$
representing critical points of the master function $
\Phi(\bs t; \bs z; \bs \La, w\cdot\tilde\La_\infty)$ 
is isomorphic to the Bruhat cell
 $\beta ( G_w^{\mc F^\infty})\ \subset \ \beta ( FL(V_P))$.
\end{cor}

\section{Selfdual spaces of polynomials}

\subsection{Dual spaces of polynomials}\label{DU}
Let $V$ be a vector space of polynomials of dimension $N + 1$ with no base
points and $D = \partial^{N+1} + \dots$ the linear differential operator with kernel $V$. 
Let $z_1, \dots , z_n,  \infty$ be  its singular points.
Let the exponents at  $z_s$ and $\infty$ be denoted
\bea
&
0,\ m_1^{(s)} + 1,\ m_1^{(s)} + m_2^{(s)} + 2,\ \dots , \ m_1^{(s)} + \dots + m_N^{(s)} + N\  ,&
\\
&
l,\ l + m_1^{(\infty)} + 1,\ l + m_1^{(\infty)} + m_2^{(\infty)} + 2,
\ \dots ,\ l + m_1^{(\infty)} + \dots + m_N^{(\infty)} + N\ , &
\eea 
respectively. Set
\bean\label{T polyn}
T_i\ =\ \prod_{s=1}^n \ (x-z_s)^{m_i^{(s)}}\ ,
\qquad i = 1, \dots , N\ .
\eean
For $u_1, \dots , u_i \in V$, define a polynomial
\be
W^{\dagger}(u_1, \dots , u_i)\ = \ W(u_1, \dots , u_i)\ \prod_{j=1}^{i-1}\  T_j^{j-i}\ ,
\ee
{\it the divided Wronskian  with respect to}
$V$. Let $\La^i (V)$ be the $i$-th exterior power of $V$. The linear map
\be
\La^i (V) \ \to \ \C[x]\ ,
\qquad
u_1 \wedge \dots \wedge u_i \ \mapsto \ W^\dagger (u_1, \dots , u_i)\ ,
\ee
is an immersion for $i\  = \  N, \ N+1$. The map  sends $\La^{N+1}(V)$ to $\C$.
Denote $V^\dagger $ the image of $\La^N(V)$. There is a non-degenerate pairing
\be
V \otimes V^\dagger  \ \to \ \C \ ,
\qquad u \otimes W^\dagger (u_1, \dots , u_N)\ 
\mapsto \ W^\dagger(u, u_1, \dots , u_{N})\ .
\ee

\begin{lem}
We have $V^{\dagger\dagger}=V$.
\end{lem}
\begin{proof} Follows from Lemma \ref{wr id 1}
\end{proof}

Let $D^\dagger = \partial^{N+1} + \dots$ be the linear differential 
operator with kernel $V^\dagger$. 
\begin{lem}
The equation  $D^\dagger u = 0$ has singular points $z_1, \dots , z_n, \infty$. The
exponents  at $z_s$ and infinity are
\bea
&
0, \ m_N^{(s)} + 1,\ m_N^{(s)} + m_{N-1}^{(s)} + 2,\ \dots ,\
m_N^{(s)} + \dots + m_1^{(s)} + N\ ,
&
\\
& l,\ l + m_N^{(\infty)} + 1,
\ l + m_N^{(\infty)} + m_{N-1}^{(\infty)} + 2,\ 
\dots ,\ l + m_N^{(\infty)} + \dots + m_1^{(\infty)} + N
\ ,
&
\eea
respectively. 
\hfill $\square$ 
\end{lem}

\begin{lem}
Let $D = D_1 \dots D_{N+1}$ be a
factorization of $D$ to linear factors, $D_i = \partial + f_i$, then 
$D^\dagger = D_{N+1} \dots D_1$.
\hfill
$\square$
\end{lem}

\subsection{Selfdual spaces and Witt bases}\label{SELF}
A space of polynomials $V$  is called {\it selfdual} if $V=V^\dagger$. 
For instance, the space of polynomials of degree not greater than $N$ is
selfdual. In this case all polynomials $T_i$ are equal to $1$.

Let $V$ be selfdual. Then 
$m_i^{(s)} = m_{N+1-i}^{(s)}$ and $m_i^{(\infty)} = m_{N+1-i}^{(\infty)}$
for all $s$ and $i$. 
There is a non-degenerate pairing 
\be
( \ ,\  )\ :\  V \otimes V \ \to\  \C\ .
\ee 
If
$u, v \in V$, then we write $v = W^{\dagger} (u_1, \dots , u_{N})$ and set 
$( u , v )\ =\ W^{\dagger}(u, u_1 , \dots , u_{N})$.
This pairing is called {\it the canonical bilinear form}.

A basis $u_1, \dots , u_{N+1}$ in a space of polynomials $V$ is called {\it
a  Witt basis} if 
\bean\label{dar 2 eqn}
u_i\ =\ W^\dagger (u_1, \dots , \widehat{u}_{N+2-i}, \dots , u_{N+1}), 
\qquad i \ =\  1, \dots , N+1 .
\eean

\begin{theorem}\label{special basis prop}
Let $V$ be selfdual. Then $V$ has a Witt basis.
\end{theorem}
\begin{cor}\label{symm}
Let $V$ be selfdual. Then the non-degenerate form $(\ ,\ )\ :\  V \otimes V \ \to\  \C$
is symmetric, if $\dim V$ is odd, and skew-symmetric, if $\dim V$ is even.
\end{cor}

{\bf Proof of Theorem \ref{special basis prop}}
Let $d_1 > d_2 > \dots >  d_{N+1}$ be degrees of polynomials in $V$,
where $d_i =  l + m_1^{(\infty)} + \dots + m_{N+1-i}^{(\infty)} + N + 1 - i$
in the notations of the previous section. Let $p_1, \dots , p_{N+1}$ be monic polynomials
in $V$, deg $p_i = d_i$. 
Denote $D_i = \prod (d_k - d_j)$ where the product is over all $k, j$ such that $k < j $ and both
$k$ and $j$ are not equal to $i$. We have $D_i = D_{N+2-i}$ since $V$ is selfdual.
\begin{lem}
For $i = 1, \dots , N+1$, we have
\be
W^\dagger (p_1, \dots , p_{i-1}, p_{i+1}, \dots , p_{N+1})\ = \
D_i \ x^{d_{N+2-i}} \ + \ \dots \ ,
\ee
where dots denote the lower degree terms.
\hfill $\square$
\end{lem}
Introduce $q_i \ = \ D_i^{-{1\over 2}}\ \prod_{j=1}^{N+1} D_j^{{1\over 2N - 2}} \ p_i$.
Then for $i = 1, \dots , N+1$, we have
\be
W^\dagger (q_1, \dots , q_{i-1}, q_{i+1}, \dots , q_{N+1})\ = \
q_{N+2-i} \ + \ \sum_{j = N+3-i}^{N+1} \ a^i_j \ q_j \ 
\ee
for suitable numbers $a_j^i$. Adding to each $q_i$  polynomials $q_{i+1}, \dots , q_{N+1}$
with suitable coefficients we can get new polynomials $q_1, \dots , q_{N+1}$ such that
\bean\label{very good}
W^\dagger (q_1, \dots , q_{i-1}, q_{i+1}, \dots , q_{N+1})\ = \
q_{N+2-i}
\eean
for $ i\ \leq \  {N\over 2} +  1$, and
\bean\label{good}
W^\dagger (q_1, \dots , q_{i-1}, q_{i+1}, \dots , q_{N+1})\ = \
q_{N+2-i} \ + \ \sum_{j = N+3-i}^{i} \ a^i_j \ q_j
\eean
for $ i\  >  \  {N\over 2} +  1$.
\begin{lem}\label{witt lem}
If the polynomials $q_1, \dots , q_{N+1}$ satisfy \Ref{very good} 
for $ i\ \leq \  {N\over 2} +  1$ and \Ref{good} for $ i\  >  \  {N\over 2} +  1$,
then they satisfy \Ref{very good} for all $i$.
\end{lem}
\begin{proof}
First let $N+1 = 2k$ be  even. 
Assume that we already showed that for some $i_0  \geq \ k $ the polynomials
$q_1, \dots , q_{N+1}$ satisfy \Ref{very good} for all $i \ < \ i_0$. 
We prove \Ref{very good} for $i = i_0$. We have
\bea\label{proof good}
\lefteqn{W^\dagger (q_1, \dots , q_{i_0-1}, q_{i_0+1}, \dots , q_{N+1})=} \\
&&{}\qquad =q_{N+2-i_0} +   a^{i_0}_{N+3-i_0}  q_{N+3-i_0}  +  \dots  + 
a^{i_0}_{i_0+1} q_{i_0+1}  + 
a^{i_0}_{i_0} q_{i_0}.
\notag
\\
\eea
We claim that $ a^{i_0}_{N+3-i_0} = \dots = a^{i_0}_{k} =
a^{i_0}_{k+1} = \dots = a^{i_0}_{i_0+1} = 0$. 

Set $W_i = W^\dagger (q_1, \dots , q_{i-1}, q_{i+1}, \dots , q_{N+1})$.
We calculate $X = W^\dagger (W_1, \dots , $ $ W_{N+1-i_0}, W_{i_0}, $
$W_{N+3-i_0}, $ $ \dots , W_k)$
in two ways. According to our assumptions we have
\be
X \ = \  W^\dagger (q_{N+1}, q_N,  \dots , q_{i_0+1},
\ q_{N+2-i_0} +  \sum_{j = N+3-i_0}^{i_0} a^{i_0}_{j}  q_{j},\
q_{i_0-1}, \dots , q_{k+1})\ .
\ee
Notice that  polynomials $W^\dagger (q_{N+1}, q_N,  \dots , q_{i_0+1},
q_{j}, q_{i_0-1}, \dots , q_{k+1})$,\ for $j = i_0$ and $j = N+2-i_0, $ $\dots ,$ $ k$, 
are non-zero, have different degrees and hence are linearly independent.

On the other hand, by Lemma \ref{wr id 1}, we have 
\be
X \ = \ \text{const}\ W^\dagger (q_{N+1}, q_N, \dots ,
 q_{i_0+1}, q_{N+2-i_0}, q_{i_0-1}, \dots , q_{k+1})\ .
\ee
Thus $a^{i_0}_{i_0}= a^{i_0}_{N+3-i_0} = \dots =  a^{i_0}_{k} = 0$.

Similarly calculating in two ways the expression
$X = W^\dagger (W_1, \dots , $ $ W_{N+1-i_0}, $ $W_{i_0}, W_{k+1},$ $ W_{k+2},
$ $ \dots ,  $ $ W_{i_0-1})$ we conclude that
$a^{i_0}_{k+1} = \dots  = a^{i_0}_{i_0-1} = 0$. Lemma \ref{witt lem} is proved for even $N+1$.

Let $N + 1 = 2k + 1$ be odd.
Assume that we already showed that for some $i_0  \geq \ k $ the polynomials
$q_1, \dots , q_{N+1}$ satisfy \Ref{very good} for all $i \ < \ i_0$. 
We prove \Ref{very good} for $i = i_0$. 

We calculate $X = W^\dagger (W_1, \dots ,  W_{k}, W_{i_0})$
in two ways. According to our assumptions we have
\be
X \ = \  W^\dagger (q_{N+1}, q_N,  \dots , q_{k+2}, 
\ q_{N+2-i_0} +  \sum_{j = N+3-i_0}^{i_0} a^{i_0}_{j}  q_{j} )\ .
\ee
The  polynomials $W^\dagger (q_{N+1}, q_N,  \dots , q_{k+1}, q_{j} )$
\ for $j = N+2-i_0, $ $\dots ,$ $ k + 1$  
are non-zero, have different degrees and hence are linearly independent.

On the other hand, we have 
\bea
X \ & = & \ \text{const}\ W^\dagger (q_{k+1}, \dots ,
q_{i_0-1}, q_{i_0+1}, \dots , q_{N+1})\ 
\\
& = &
W^\dagger (W_{k-1}, \dots ,
W_{N+3-i_0}, W_{N+1-i_0}, \dots , W_1)
\\
& = & \text{const} \ W^\dagger (q_{N+2-i_0}, q_{k+1}, \dots , q_{N+1})\ . 
\eea
Thus $ a^{i_0}_{N+3-i_0} = \dots =  a^{i_0}_{k+1} = 0$.

Similarly calculating in two ways the expression
\be
X\  =\  W^\dagger (W_1, W_2, \dots , W_{N+2-i_0}, W_{k+2}, W_{k+3}, \dots , W_{i_0})
\ee
we show that $ a^{i_0}_{k+2} = \dots =  a^{i_0}_{i_0} = 0$.
Lemma \ref{witt lem} is proved for odd $N+1$.
\end{proof}
Theorem \ref{special basis prop} is proved.
\hfill
$\square$

\subsection{Isotropic flags}\label{isotropic}
Let $u_1, \dots , u_{N+1}$ be a basis in a vector space
$V$ of polynomials. Denote $W_i\ =\ W^\dagger (u_1, \dots , \widehat{u}_{i}, 
\dots , u_{N+1})$. 
\begin{theorem}\label{suf self}
Assume that
\bean\label{dar 1 eqn}
W^\dagger (u_1, \dots , u_i)\ =\ a_i\ W^\dagger (u_1, \dots , u_{N+1-i}),
\qquad i\ =\ 1, \dots , N , 
\eean
for some non-zero numbers $a_1, \dots , a_N$. Then for every $i$, the polynomial
$u_i$ is a linear combination of $W_{N+1}, W_N, \dots , W_{N+2-i}$.
\end{theorem}
\begin{proof} 
The proof is  by induction on $i$. 
For $i=1$, the polynomial $u_1$ is proportional to $W_{N+1}$ by assumption.
Let the theorem be proved for $i = 1, \dots , i_0 - 1$. Then
using Lemma \ref{wr id 1} we obtain
\bea
W^\dagger (u_1, \dots , u_{i_0-1}, W_{N+2-i_0})\ &=&\
\text{const}\  W^\dagger (W_{N+1}, W_{N}, \dots ,  W_{N+2-i_0})
\\
&=&\ \text{const}\ 
W^\dagger (u_1, \dots , u_{N+1-i_0})\ (W^\dagger (u_1, \dots , u_{N+1}))^{i_0-1} .
\eea
Notice that $W^\dagger (u_1, \dots , u_{N+1})$ is a constant.
Therefore using \Ref{dar 1 eqn} we obtain
\be
  W^\dagger (u_1, \dots , u_{i_0-1}, W_{N+2-i_0})
\ =\ \text{const} \ W^\dagger(u_1, \dots , u_{i_0}) . 
\ee
Hence $W_{N+2-i_0}$ is a linear combination of $u_1, \dots , u_{i_0}$ and the theorem is proved.
\end{proof}
\begin{cor} \label{self cor}
If $V$ has a basis satisfying \Ref{dar 1 eqn}, then $V$ is selfdual.
\end{cor}
Let $V$ be a selfdual space of polynomials. 
For a subspace $U \subset V$ denote $U^\perp$ its  orthogonal
complement,  $( U , U^\perp ) = 0$,
 $\dim U  + \dim U^\perp = \dim V $.
A full flag $\mathcal F = \{F_1 \subset \dots \subset F_{N+1} = V\}$ 
is called {\it isotropic}  if $F_i^{\perp} = F_{N+1-i}$ for \ $i = 1, \dots , N$. 

\begin{cor}
Let $u_1, \dots , u_{N+1}$ be a basis in $V$ satisfying \Ref{dar 1 eqn}, then the 
full flag generated by the basis is isotropic.
\end{cor}

\begin{lem}
Let $\mc F$ be an isotropic flag in $V$ and $u_1, \dots , u_{N+1}$  a basis in $V$ 
adjusted to $\mc F$. Then the basis satisfies \Ref{dar 1 eqn}. 
\end{lem}
\begin{proof}
Since $F_i^\perp = F_{N+1-i}$, we have two bases in $F_{N+1-i}$: 
the basis $u_1, \dots , u_{N+1-i}$ and the basis
$W_{N+1}, W_N, \dots , W_{i+1}$. Hence
\bea
W^\dagger (u_1, \dots , u_{N+1-i}) \ = 
\ \text{const}\ W^\dagger (W_{N+1}, W_N, \dots , W_{i+1})\
=\ \ \text{const}\ W^\dagger (u_1,  \dots , u_i)\ .
\eea
\end{proof}

\subsection{Generating procedures for isotropic flags}\label{gen iso}
Let $V$ be a selfdual space of polynomials. 
Let $FL(V)$ be the variety of all full flags of $V$. Denote $FL^\perp(V) \subset FL(V)$ the 
subvariety of all isotropic flags.

Let  $\beta : FL(V) \to  \PCN ,\ \mathcal F  \mapsto  y^{\mathcal F}$, be
the generating isomorphism associated with $V$,
see Section \ref{ha-ha}. Remind that if $u_1, \dots , u_{N+1}$ is a basis adjusted to
a flag $\mc F \in FL(V)$, then $\bs y^{\mc F} = (y_1^{\mc F}, \dots , y_N^{\mc F})$, where
$y_i^{\mc F} = W^\dagger (u_1, \dots , u_i)$.
According to Section \ref{isotropic}, a flag $\mc F$ is isotropic if and only its 
image under the generating isomorphism is symmetric:
\bean\label{symm odd}
y^{\mc F}_i \ =\ y^{\mc F}_{N+1-i} \qquad \text{for} 
\qquad
i = 1, \dots , N\ .
\eean
The group of linear transformations of $V$ preserving the scalar product $(\ ,\ )$ acts on the
space of isotropic flags $FL^\perp (V)$. We describe 
the infinitesimal action of this group on the image $\beta (FL^\perp(V))$.

First we assume that $N+1 = 2k$ is even. In this case the scalar product 
$(\ , \ )$ is skew-symmetric. 
The special symplectic Lie algebra of $V$ consists of all
traceless endomorphisms $x$ of $V$
such that $(xv , v') + (v , xv') = 0$ for all $v, v' \in V$.

Let $\bs u = (u_1, \dots , u_{2k})$ be a Witt basis in $V$. We have
\be
( u_i \ , \ u_{2k+1-i})\ =\ (-1)^{i+1} , \qquad i = 1, \dots , k ,
\ee
and $(u_i , u_j) = 0$ if $i + j \neq 2k+1$.
The choice of the basis identifies the special symplectic Lie algebra 
with a Lie subalgebra of $sl_{2k}$, which is denoted
$sp_{2k}$. The Lie algebra $sp_{2k}$ has the root system of type $C_k$.

Denote  $E_{i, j}$  the matrix with zero entries except 1 at the intersection of
the $i$-th row and $j$-th column.

The lower triangular part of $sp_{2k}$ is spanned by
matrices $ \ E_{i, i+1}\ + \ E_{ 2k-i, 2k+1-i}$  for $i = 1, \dots, k - 1$
and $E_{k, k+1}$. Denote these matrices $X_1, \dots , X_k$, respectively.
The matrices define linear transformations of $V$.

For any $i \in \{1, \dots , k\}$ and $c \in \C$, the basis $e^{c X_i} \bs u$ is a Witt basis.
Let $\mc F$ be the isotropic flag generated by $\bs u$ and $e^{c X_i} \mc F$  the 
isotropic flag generated by $e^{c X_i} \bs u$. We describe the dependence on $c$ of the
$2k-1$-tuple $\beta(e^{c X_i} \mc F)$.

For $i < k$, we have
\be
e^{c X_i} \bs u \ = \
(u_1, \dots , u_{i-1}, u_i + c u_{i+1}, u_{i+1}, \dots ,
u_{2k-i}, u_{2k+1-i} + c u_{2k+2-i}, u_{2k+2-i}, \dots , u_{2k})
\ee
and
\bea
\lefteqn{\beta (e^{c X_i} \mc F) = }\\
&&=(y^{\mc F}_1, \dots, y^{\mc F}_{i-1},  y_i(x, c),  y^{\mc F}_{i+1},
\dots , y^{\mc F}_{k-1}, y^{\mc F}_k, y^{\mc F}_{k-1}, \dots , y^{\mc F}_{i+1}, y_i(x, c),  y^{\mc F}_{i-1}, \dots ,  y^{\mc F}_{1}) ,
\eea
with $y_i(x, c) = W^\dagger (u_1, \dots , u_{i-1},  u_i + c u_{i+1})$.
\begin{lem}\label{skew 1} 
For $i < k$, we  have
\be
W (y_i(x, c), {\partial y_i \over \partial c}(x, c))\ =\
T_i\ y^{\mc F}_{i-1}\ y^{\mc F}_{i+1}\ .
\ee
\hfill
$\square$
\end{lem}
For $i = k$, we have
\be
e^{c X_i} \bs u \ = \
(u_1, \ \dots ,\ u_{k-1},\ u_k + c u_{k+1},\ u_{k+1},
\  \dots ,\ u_{2k})
\ee
and
\be
\beta(e^{c X_i} \mc F) \ = \ (y^{\mc F}_1, \dots,\ y^{\mc F}_{k-1},\ y_k(x, c), 
 \ y^{\mc F}_{k-1}, \ \dots , y^{\mc F}_1)
\ee
with $y_k(x, c)\ = \ W^\dagger (u_1, \dots , u_{k-1},  
u_k + c u_{k+1})$.
\begin{lem} \label{skew 2} 
We have
\be
W (y_k(x, c), {\partial y_k \over \partial c}(x, c))\ =\
T_k \ (y^{\mc F}_{k-1})^2  \ .
\ee
\hfill
$\square$
\end{lem}
Notice that in Lemmas \ref{skew 1} and \ref{skew 2} the function
$ {\partial y_k \over \partial c}(x, c)$ does not depend on $c$.

Now we assume that $N+1 = 2k+1$ is odd. In this case the scalar product 
$(\ , \ )$ is symmetric. The special orthogonal Lie algebra of $V$ consists of all
traceless endomorphisms $x$ of $V$
such that $(xv , v') + (v , xv') = 0$ for all $v, v' \in V$.

Let $\bs u = (u_1, \dots , u_{2k+1})$ be a Witt basis in $V$. We have
\be
( u_i \ , \ u_{2k+2-i})\ =\ (-1)^{i+1} , \qquad i = 1, \dots , k+1 ,
\ee
and $(u_i , u_j) = 0$ if $i + j \neq 2k+2$.
The choice of the basis identifies the special orthogonal Lie algebra 
with a Lie subalgebra of $sl_{2k+1}$, which is  denoted
$so_{2k+1}$. The Lie algebra $so_{2k+1}$ has the root system of type $B_k$.

The lower triangular part of $so_{2k+1}$ is spanned by
matrices $E_{i, i+1}\ + \ E_{ 2k+1-i, 2k+2-i}$  for $i = 1, \dots, k-1$
and $E_{k, k+1} \ +\ E_{ k+1, k+2}$. Denote these matrices $X_1, \dots , X_k$, respectively.
The matrices define linear transformations of $V$.

For any $i \in \{1, \dots , k\}$ and $c \in \C$, the basis $e^{c X_i} \bs u$ is a Witt basis.
Let $\mc F$ be the isotropic flag generated by $\bs u$ and $e^{c X_i} \mc F$  the 
isotropic flag generated by $e^{c X_i} \bs u$. We describe the dependence on $c$ of the
$2k$-tuple $\beta(e^{c X_i} \mc F)$.

For $i < k$, we have
\bea
\lefteqn{e^{c X_i} \bs u \ =}\\
&&(u_1, \dots , u_{i-1}, u_i + c u_{i+1}, u_{i+1}, \dots ,
u_{2k+1-i}, u_{2k+2-i} + c u_{2k+3-i}, u_{2k+3-i}, \dots , u_{2k+1})
\eea
and
\bea
\lefteqn{\beta (e^{c X_i} \mc F)= }\\
&&=(y^{\mc F}_1,\dots, y^{\mc F}_{i-1}, 
y_i(x, c), y^{\mc F}_{i+1}, 
 \dots ,y^{\mc F}_k, y^{\mc F}_k,\dots ,y^{\mc F}_{i+1},y_i(x, c),
 y^{\mc F}_{i-1}, \dots ,  y^{\mc F}_{1}) ,
\eea
with $y_i(x, c) = W^\dagger (u_1, \dots , u_{i-1},  u_i + c u_{i+1})$.
\begin{lem}\label{sym 1}  
For $i < k$ we have
\be
W (y_i(x, c), {\partial y_i \over \partial c}(x, c))\ =\
T_i\ y^{\mc F}_{i-1}\ y^{\mc F}_{i+1}\ .
\ee
\hfill
$\square$
\end{lem}
For $i = k$, we have
\be
e^{c X_i} \bs u \ = \
(u_1, \dots , u_{k-1},\ u_k + c u_{k+1} + {c^2\over 2} u_{k+2},\ u_{k+1} + c u_{k+2},
\ u_{k+2}, \dots , u_{2k+1})
\ee
and
\be
\beta(e^{c X_i} \mc F) \ = \ (y^{\mc F}_1,\ \dots,\ y^{\mc F}_{k-1},\ y_k(x, c), 
\ y_k(x, c), \ y^{\mc F}_{k-1}, \ \dots ,\ y^{\mc F}_1)
\ee
with $y_k(x, c)\ = \ W^\dagger (u_1, \dots , u_{k-1},  
u_k + c u_{k+1} + {c^2\over 2} u_{k+2})$.
\begin{lem}\label{sym 2}  
We have
\bean\label{important}
W (y_k(x, c), {\partial y_k \over \partial c}(x, c))\ =\
T_k \ y^{\mc F}_{k-1}\ y_k(x, c) \ .
\eean
\hfill
$\square$
\end{lem}
\begin{lem}\label{square}
Assume that the two polynomials $y^{\mc F}_{k}$ and $T_k \,y^{\mc F}_{k-1}$
do not have common roots. Then $y^{\mc F}_{k} = p^2$,
$y_k(x, c) = (p + c q)^2$  for  suitable polynomials $p(x), \, q(x)$ and 
\bean\label{also important}
W (p\ , \ q )\ =\
2 \ T_k \ y^{\mc F}_{k-1}\  .
\eean
\end{lem}
\begin{proof}
The polynomial $y_k(x, c)$ is a polynomial in $x$ and $ c$, quadratic in $c$.
We claim that $y_k$ cannot be of the form $p(x) q(x,c)$ where $p$ is a polynomial 
in $x$ of positive degree and $q$ is a polynomial in $x, c$. Indeed, in this
case equation \Ref{important} would imply that the roots of $p$ lie among the roots of 
$T_k \,y^{\mc F}_{k-1}$.

We claim that $y_k$ cannot be irreducible. Indeed,  then there 
would exist $(z_0, c_0) \in \C^2 $ and a
germ of a holomorphic function $ f : (\C , c_0) \to
(\C, z_0)$ such that in a neighborhood of $(z_0, c_0)$ the polynomial $y_{k}$ would be
presentable in the form $y_k = (x - f(c)) g(x, c)$ with $g(z_0, c_0) \neq 0$.
But this presentation would contradict to \Ref{important}, because its left hand side and 
right hand side would be having different orders of zero with respect to $x$ at $x = f(c)$.

The same reason shows that $y_k$ cannot be presented in the form
$p(x, c) q(x,c)$ with  polynomials $p, q$ not proportional.
Thus $y_k = (p(x) + c q(x))^2$ for suitable polynomials $p$ and $q$. Substituting 
this expression to \Ref{important} we get \Ref{also important}.
\end{proof}

\section{The case of root systems of type $B_N$ and $C_N$ }

\subsection{Critical points of type $B_N$}\label{BN}
Consider the root system of type $B_N$. Let
$\al_1, \dots , \al_{N-1}$ be the long simple roots and $\al_N$ the short one.
 We have 
\bea
(\al_N , \al_N)\ =\ 2 , \qquad
(\al_i , \al_i)\ =\ 4 , \qquad
(\al_i , \al_{i+1}) \ =\ - 2 , \qquad i = 1, \dots , N-1 ,
\eea
and all other scalar products are equal to zero.
The root system $B_N$ corresponds to the Lie algebra $so_{2N+1}$. We denote $\h$ 
its Cartan subalgebra.

We consider also the root system of type $A_{2N-1}$ with simple roots $\al^A_1, \dots ,
\al^A_{2N-1}$. The root system $A_{2N-1}$ corresponds to the Lie algebra $sl_{2N}$. We
denote $\h_A$ its Cartan subalgebra.

We have a map $\h^* \to \h_A^*, \  \La \mapsto \La^A$, where $\La^A$ is defined by
\be
\langle \La^A\  , \ (\al_i^A)^\vee \rangle\ =\ \langle \La^A\ ,\ (\al_{2N-i}^A)^\vee
\rangle \ =\
\langle \La \ , \ (\al_i)^\vee \rangle , \qquad i = 1 , \dots ,  N . 
\ee

Let $\La_1, \dots, \La_n \in \h^*$ be dominant integral $so_{2N+1}$-weights,
$z_1, \dots , z_n$ complex numbers. Let the polynomials $T_1, \dots , T_N$ be given by \Ref{T}.
Remind that an $N$-tuple
  of polynomials $\bs y$ represents a critical point of a master function
associated with $so_{2N+1}$, \ $\La_1, \dots , \La_n$, \ $z_1, \dots , z_n$,
if and only if $\bs y$ is generic with respect to weights
\ $\La_1, \dots , \La_n$ of $so_{2N+1}$, and points $z_1, \dots z_n$
and
there exist polynomials $\tilde y_i$,\  $i = 1, \dots , N$,  such that
\bea
&& 
W(y_i \ ,\  \tilde y_i )\ =\ T_i\ y_{i-1}\ y_{i+1} ,\qquad i\ =\ 1, \dots , N-1 ,
\\
&&
W(y_N \ , \ \tilde y_N)\ =\ T_N\ y_{N-1}^2 .
\eea

For an $N$-tuple of polynomials $\bs y = (y_1, \dots, y_N)$, let 
$\bs y_A$ be the $2N-1$-tuple of polynomials
$(y_1, \dots , y_{N-1}, y_N, y_{N-1}, \dots, y_1)$.

\begin{lem}
An $N$-tuple $\bs y$ represents a critical point of a master function
associated with  $so_{2N+1}$, \ $\La_1, \dots , \La_n$, \ $z_1, \dots , z_n$,
if and only if the $2N-1$-tuple of polynomials
$\bs y_A$ represents a critical point of a master function
associated with $sl_{2N}$,\ $\La_1^A, \dots ,$ $ \La_n^A$, \ $z_1, $ $\dots ,$ $ z_n$.
\hfill $\square$
\end{lem}

Let $\bs y^0$ represent a critical point of a master function
associated with  $so_{2N+1}$, \ $\La_1, \dots , \La_n$, \ $z_1, \dots , z_n$.
Denote $P$ the $so_{2N+1}$-population of critical points originated at $\bs y^0$. 
Denote $P^A$ the $sl_{2N}$-population of critical points originated at $\bs y^0_A$. 
Clearly there is an injective map $P \to P^A$, $\bs y \mapsto \bs y_A$.

\begin{theorem}\label{so-self B}
The fundamental space of $P^A$ is selfdual.
\end{theorem} 
\begin{proof}
Denote $V$ the fundamental space of $P^A$.
Let $u_1, \dots , u_{2N}$ be a basis in $V$ 
constructed using $\bs y_A$ as in Section \ref{all polynom section}.
By Lemma \ref{wr(u) lem} we have
\bea
W^\dagger (u_1, \dots , u_i)\ =\ W^\dagger (u_1, \dots , u_{2N-i}),
\qquad i\ =\ 0, \dots , 2N .
\eea
The theorem follows from Corollary \ref{self cor}.
\end{proof}
The fundamental space of $P^A$ is called {\it the fundamental space} of the 
$so_{2N+1}$-population $P$ and denoted $V_P$. We have $\dim V_P = 2N$.

\begin{theorem}\label{cor so}
The generating morphism $\beta : FL(V_P) \to \bs P (\C[x])^{2N-1}$ identifies
the subvariety $FL^\perp(V_P)$ of flags in $V_P$, isotropic with respect to the skew-symmetric
form $(\ ,\ )$, and the $so_{2N+1}$-population $P$.
\end{theorem}
\begin{proof}
Let $\bs y=\beta(\mc F)$ be in $P$, where $\mc F$ is some isotropic
flag. Then the $N$-tuples 
$\beta (e^{cX_i}\mc F)=(y_1(x),\dots, y_i(x,c),\dots,y_N(x))$ are 
the immediate descendants of $\bs y$ in the $i$-th direction. Indeed,
this follows from  Lemmas \ref{skew 1} and \ref{skew 2} with $c=0$ since we 
have $y_i(x,c)=y_i+c\partial y_i/\partial c(x,0)$. 

Therefore the generating procedures
for the $so_{2N+1}$-population $P$ coincide with the generating procedure for the
variety of full flags in the fundamental space $V_P$ isotropic with respect to the skew-symmetric
form  $(\ ,\ )$. 
\end{proof}

The Lie algebra $so_{2N+1}$ is Langlands dual to the Lie algebra $sp_{2N}$. According 
to the corollary, a population of critical points associated with $so_{2N+1}$ is 
isomorphic to the full flag variety of its Langlands dual $sp_{2N}$.

Let $P$ be an $so_{2N+1}$-population associated with integral dominant weights
$\La_1, \dots \La_n$ and points $z_1, \dots , z_n$. Let $\bs y \in P$. Introduce a linear 
differential operator $D(\bs y)$ of order $2N$ by the formula
\bea
D(\bs y)\ &=&\
(\partial - \ln'({T_1^2 \dots T_{N-1}^2 T_N\over y_1}))
(\partial - \ln'({y_1 T_1^2 \dots T_{N-1}^2 T_N\over y_2 T_1}))
\\
&&
(\partial - \ln'({y_2 T_1^2 \dots T_{N-1}^2 T_N\over y_3 T_1 T_2}))
\ \dots\
(\partial - \ln'({y_{N-1} T_1\dots T_{N-1} T_N\over y_N}))
\\
&&
(\partial - \ln'({y_{N} T_1\dots T_{N-1}\over y_{N-1}}))
(\partial - \ln'({y_{N-1} T_1\dots T_{N-2}\over y_{N-2}}))
\dots
(\partial - \ln'(y_{1}))
\eea
where $T_1, \dots T_N$ are as above.
\begin{cor}
The operator $D(\bs y)$ does not depend on the choice of $\bs y$ in $P$. Its kernel 
 is the fundamental space $V_P$.
\end{cor}
The operator is called {\it associated with the population} and is denoted $D_P$.

By Theorem \ref{cor so} any $so_{2N+1}$-population is isomorphic to the variety of isotropic
flags of the fundamental space of the population. The converse is also true.

\begin{theorem}\label{conv so}
Let $V$ be a selfdual space of polynomials of even dimension $2N$. Let $FL^\perp(V)$ be the variety
of the isotropic full flags of $V$, and
$\beta (FL^\perp(V))$ its image under the generating morphism. Then
$\beta (FL^\perp(V))$ is an $so_{2N+1}$-population with fundamental space $V$.
\end{theorem}

Theorem \ref{conv so} is proved in Section \ref{sec loc}.

\subsection{Critical points of type $C_N$}\label{CN}
Consider the root system of type $C_N$. Let
$\al_1, \dots , \al_{N-1}$ be the short simple roots and $\al_N$ the long one.
We have 
\bea
(\al_N , \al_N)\ =\ 4 ,  \qquad && \qquad (\al_{N-1} , \al_N)\ =\ -\ 2  , 
\\
(\al_i , \al_i)\ =\ 2 , \qquad
 (\al_i , \al_{i+1}) \ &=&\ - \ 1 ,  \qquad i = 1, \dots , N-2 ,
\eea
and all other scalar products are equal to zero.
The root system $C_N$ corresponds to the Lie algebra $sp_{2N}$. We denote $\h$ 
its Cartan subalgebra.

We consider also the root system of type $A_{2N}$ with simple roots $\al^A_1, \dots ,
\al^A_{2N}$. The root system $A_{2N}$ corresponds to the Lie algebra $sl_{2N+1}$. We
denote $\h_A$ its Cartan subalgebra.

We have a map $\h^* \to \h_A^*, \  \La \mapsto \La^A$, where $\La^A$ is defined by
\be
\langle \La^A\  , \ (\al_i^A)^\vee \rangle\ =\ \langle \La^A\ ,\ (\al_{2N+1-i}^A)^\vee
\rangle \ =\
\langle \La \ , \ (\al_i)^\vee \rangle , \qquad i = 1 , \dots ,  N . 
\ee

Let $\La_1, \dots, \La_n \in \h^*$ be dominant integral $sp_{2N}$-weights,
$z_1, \dots , z_n$ complex numbers. Let the polynomials 
$T_1, \dots , T_N$ be given by \Ref{T}.
Remind that an $N$-tuple of polynomials $\bs y$ represents a critical 
point of a master function associated with $sp_{2N}$,
 \ $\La_1, \dots , \La_n$, \ $z_1, \dots , z_n$,
if and only if $\bs y$ is generic with respect to
weights  \ $\La_1, \dots , \La_n$ of $sp_{2N}$ and points \ $z_1, \dots , z_n$,
and there exist polynomials $\tilde y_i$, \ $i = 1, \dots , N$,  such that
\bea
W(y_i \ ,\  \tilde y_i )\ =\ T_i\ y_{i-1}\ y_{i+1} ,&\qquad &
i\ =\ 1, \dots , N-2 , 
\\
W(y_{N-1} \ ,\ \tilde y_{N-1})\ =\ T_{N-1}\ y_{N-2}\ y_{N}^2 ,
&
\qquad
&
W(y_{N} \ ,\  \tilde y_{N})\ =\ T_{N}\ y_{N-1} .
\eea

For an $N$-tuple of polynomials $\bs y = (y_1, \dots, y_N)$, let 
$\bs y_A$ be the $2N$-tuple of polynomials
\be
\bs y_A\ =\ (y_1, \dots ,\ y_{N-1},\  y^2_{N},\  y^2_N,\  y_{N-1}, \ \dots, y_1) .
\ee

Let $\bs y$ represent a critical 
point of a master function associated with $sp_{2N}$,
 \ $\La_1, \dots , \La_n$, \ $z_1, \dots , z_n$.
We shall construct an $sl_{2N+1}$-population containing $\bs y_A$.
Set 
\bea
\bs y_{A,1}\ & = &\
(y_1, \ \dots ,\ y_{N-1},\  y_{N} \tilde y_N,\  y^2_N,\  y_{N-1}, \dots,\ y_1)\ ,
\\
\bs y_{A,2}\ & = &\
(y_1, \ \dots ,\ y_{N-1},\  y_{N} \tilde y_N,\ y_N^2 + c (\tilde y_N)^2,\  y_{N-1}, \dots,\ y_1)\ ,
\eea
where $c$ is a non-zero number.
\begin{lem}\label{sp pop}
For almost all $c$, the tuple $\bs y_{A,2}$ represents a critical point of a master 
function associated with $sl_{2N+1}$, $\La^A_1, \dots , \La_n^A$, $z_1, \dots , z_n$.
The $sl_{2N+1}$-population,
containing the critical point $\bs y_{A, 2}$, does not depend
on $c$ and contains also the tuple $\bs y_A$. 
\end{lem}
\begin{proof}
Since $\bs y$ is generic with respect to weights $\La_1, 
\dots , \La_n$ of $sp_{2N}$, and points $z_1, \dots , z_n$, clearly
$\bs y_{A,2}$ is generic with respect to weights
$\La^A_1, \dots , \La_n^A$ of $sl_{2N+1}$, and points $z_1, \dots , z_n$
for almost all $c$.
We have 
\be
W( y^2_N \ ,\ y_{N} \tilde y_N )\ =\ T_N \ y_{N-1}\ (y_N)^2 ,
\qquad
W( (y_N)^2\ , \ y^2_N + c(\tilde y_N)^2)\ =\ 2\ c\ T_N \ y_{N-1}\  \ y_{N} \tilde y_N .
\ee
The first equation shows that $\bs y_A$ is fertile with respect to
$sl_{2N+1}$, $\La^A_1, \dots , \La_n^A$, $z_1, \dots , z_n$.
The equations also show that with respect to
$sl_{2N+1}$, $\La^A_1, \dots , \La_n^A$, $z_1, \dots , z_n$ the tuple
 $\bs y_{A,1}$ is fertile in all directions but the $N-1$-th. Now we  show that
$\bs y_{A,1}$ is fertile in direction $N-1$.

Let $y_{N-1} = \prod_i (x-t_i),\ y_{N} = \prod_i (x-u_i),\
 \tilde y_{N} = \prod_i\ (x-v_i)$.
The numbers $\{t_i\}$ form an $sl_2$ critical point with weight 
$y_{N-2} T_{N-1} y_{N}^2$.
Therefore we have its critical point  equations 
\be
\sum_j\ {2\over t_i - u_j }\ +\ X\ = \ 0,
\ee
where $X$ denotes terms which do not depend on $y_N$.
Similarly, 
$\{t_i\}$ form an $sl_2$ critical point with weight
$y_{N-2} T_{N-1} \tilde y_{N}^2$.
Therefore we have its critical point equations
\be
\sum_j\ {2\over t_i - v_j}\ +\ X\ =\ 0
\ee 
with the same $X$. Hence
\be
\sum_j\ {1\over t_i - u_j }\ +\ \sum_j\ {1\over t_i - v_j}\ +\ X\ =\ 0. 
\ee
This means that $\{t_i\}$ form an $sl_2$  critical point with weight
$y_{N-2} T_{N-1} y_N \tilde y_N$. Thus $\bs y_{A,1}$ is fertile in direction $N-1$.

With respect to 
$sl_{2N+1}$, $\La^A_1, \dots , \La_n^A$, $z_1, \dots , z_n$ the tuple
$\bs y_{A,2}$ is fertile in all directions but the $N$-th
and $N+2$-th. We show that
$\bs y_{A,2}$ is fertile in these directions too.

We have
\be
W(y_{N-1} , \bar y_{N-1})\ =\ T_{N-1}\ y_{N-2}\ y_N^2\ ,
\qquad
W(y_{N-1} , \tilde y_{N-1})\ =\ T_{N-1}\ y_{N-2}\ \tilde y_N^2
\ee
for suitable polynomials $\bar y_{N-1}$ and $\tilde y_{N-1}$. Adding the equations we get
$W(y_{N-1}\ , \ \bar y_{N-1} + c \tilde y_{N-1})\ =\ T_{N-1}\ y_{N-2}\ ( y_N^2 + c
\tilde y^2_N)$ which means fertility in direction $N+2$. Similarly, we have
$W( y_{N}^2 + c \tilde y^2_{N}\ , \ y_{N} \tilde y_N )\ =\ T_{N}\ y_{N-1}\ ( y_N^2 + c
\tilde y^2_N)$ which means fertility in direction $N$. These reasons show that
for all $c$ the tuples $ \bs y_{A}$
$ \bs y_{A, 1}, \bs y_{A, 2}$ belong to the $sl_{2N+1}$-population originated at $\bs y_{A,2}$.
\end{proof}
Let $\bs y$ represent a critical 
point of a master function associated with $sp_{2N}$,
 \ $\La_1, \dots , \La_n$, \ $z_1, \dots , z_n$.
Let $P$ be the $sp_{2N}$-population originated at $\bs y$. Let $P^A$ be the
$sl_{2N+1}$-population, constructed in Lemma \ref{sp pop} and containing $\bs y_A$.

\begin{theorem}\label{so-self C}
The fundamental space of $P^A$ is selfdual.
\end{theorem} 
\begin{proof}
Denote $V$ the fundamental space of $P^A$.
Let $u_1, \dots , u_{2N+1}$ be a basis in $V$ 
constructed using $\bs y_A$ as in Section \ref{all polynom section}.
By Lemma \ref{wr(u) lem} we have
\bea
W^\dagger (u_1, \dots , u_i)\ =\ W^\dagger (u_1, \dots , u_{2N-i}),
\qquad i\ =\ 0, \dots , 2N + 1 .
\eea
The theorem follows from Corollary \ref{self cor}.
\end{proof}

The fundamental space of $P^A$ is called {\it the fundamental space} of the 
$sp_{2N}$-population $P$ and denoted $V_P$. We have $\dim V_P = 2N + 1$.

\begin{theorem}\label{cor sp}
The generating morphism $\beta : FL(V_P) \to \bs P (\C[x])^{2N}$ identifies
the subvariety $FL^\perp(V_P)$ of flags in $V_P$, isotropic with respect to the symmetric
form $(\ ,\ )$, and the $sp_{2N}$-population $P$.
\end{theorem}
\begin{proof}
Follows from  Lemmas \ref{sym 1},  \ref{sym 2} with $c=0$ and
\ref{square}, cf. proof of Theorem \ref{cor sp}. 
\end{proof}

The Lie algebra $sp_{2N}$ is Langlands dual to the Lie algebra $so_{2N+1}$. According 
to the corollary, a population of critical points associated with $sp_{2N}$ is 
isomorphic to the full flag variety of its Langlands dual $so_{2N+1}$.

Let $P$ be an $sp_{2N}$-population associated with integral dominant weights
$\La_1, \dots \La_n$ and points $z_1, \dots , z_n$. Let $\bs y \in P$. Introduce a linear 
differential operator $D(\bs y)$ of order $2N+1$ by the formula
\bea
D(\bs y)\ &=&\
(\partial - \ln'({T_1^2 \dots T_{N}^2\over y_1}))
(\partial - \ln'({y_1 T_1^2 \dots T_{N}^2\over y_2 T_1}))
\dots
(\partial - \ln'({y_{N-2} T_1^2 \dots T_{N}^2\over y_{N-1} T_1 \dots T_{N-2}}))
\\
&&
(\partial - \ln'({y_{N-1} T_1^2 \dots T_{N}^2\over y_N^2
T_1 \dots T_{N-1}}))
(\partial - \ln'(T_1\dots T_{N}))
(\partial - \ln'({y_{N}^2 T_1\dots T_{N-1}\over y_{N-1}}))
\\
&&(\partial - \ln'({y_{N-1} T_1\dots T_{N-2}\over y_{N-2}}))
\dots
(\partial - \ln'({y_{2} T_1\over y_{1}}))
(\partial - \ln'(y_{1}))
\eea
where $T_1, \dots T_N$ are as above.
\begin{cor}
The operator $D(\bs y)$ does not depend on the choice of $\bs y$ in $P$. Its kernel 
 is the fundamental space $V_P$.
\end{cor}
The operator is called {\it associated with the population} and is denoted $D_P$.

According to Theorem \ref{cor sp} any $sp_{2N}$-population is isomorphic to the variety of isotropic
flags of the fundamental space of the population. The converse is also true.

\begin{theorem}\label{conv sp}
Let $V$ be a selfdual space of polynomials of odd dimension $2N+1$. Let $FL^\perp(V)$ be the variety
of the isotropic full flags of $V$, and $\beta (FL^\perp(V))$ its image under the generating morphism. Then
$\beta (FL^\perp(V))$ is an $sp_{2N}$-population with fundamental space $V$.
\end{theorem}

Theorem \ref{conv sp} is proved in Section \ref{sec loc}.

\subsection{The $so_{2N+1}$-populations and Bruhat cells}
Let $P$ be an $so_{2N+1}$-population with fundamental space $V_P$.
The space $V_P $ is a  selfdual space of dimension $2N$.
The canonical bilinear form of $V_P$ is skew-symmetric. Let $z_1, \dots , z_n,
\infty$ 
be ramification points of $V_P$.
Let $\La_1^A, \dots , \La_n^A, \tilde\La_\infty^A$
be the corresponding integral dominant $sl_{2N}$-weights of $V_P$ and 
$\La_1^B, \dots  , \La_n^B, \tilde \La_\infty^B$ the corresponding integral dominant 
$so_{2N+1}$-weights of $P$.

The generating morphism $\beta$ identifies the variety of isotropic flags 
$FL^\perp (V_P)$ with the population $P$.
The variety $FL^\perp(V_P)$ is isomorphic to the flag variety of the Lie group 
$SP_{2N}$. 
The root system of its Lie algebra $sp_{2N}$ is of type $C_N$.

For any semisimple Lie group $G$ with a Borel subgroup $B$,
the flag variety $G/B$ has a Bruhat cell decomposition
$G/B\  = \ \sqcup_{w\in W}\, BwB/B$, see \cite{H}. 

A Bruhat cell decomposition of the flag variety of
$SL_{2N}$ was described in Section \ref{Bruhat section} 
as a cell decomposition of the variety of full flags $FL(V_P)$ with respect to the
flag ${\mc F}^\infty$. 
We have
\be
FL(V_P)\ =\ \sqcup_{w\in S^{2N}}\ G_w^{{\mc F}^\infty}.
\ee
Now we describe a Bruhat cell decomposition of the flag variety of  $SP_{2N}$
as a decomposition of $FL^\perp(V_P)$.

Let $\iota: \h^*_{sp_{2N}}\to \h^*_{sl_{2N}}$ be the linear inclusion 
given by
\be
\sqrt{2} \al_i^C \ \mapsto  \ \al_i^A + \al_{2N-i}^A, \qquad \sqrt{2} \al_N^C\ \mapsto\ 2 \al_N^A,
\ee
$i = 1, \dots , N-1$.  We have
$(\al_i^C , \al_j^C) = (\iota(\al_i^C) , \iota(\al_j^C))$ and therefore 
$\iota$ is orthogonal.
We have the following relation between simple reflections:
\be
\iota (s_i^C\La)\ =\ s_i^A s^A_{2N-i} \iota(\La) ,
\qquad 
\iota (s_N^C\La)\ =\ s_N^A \iota(\La),
\ee
$i = 1, \dots , N-1$, for any weight $\La \in \h^*_{sp_{2N}}$.
Therefore the map from the $sp_{2N}$ Weyl group to the $sl_{2N}$ Weyl group
given by 
\be
s_i^C\mapsto s^A_is^A_{2N-i},\qquad s_N\mapsto s_N,
\ee
$i=1,\dots,N-1$, is a group homomorphism. The map identifies the
$sp_{2N}$ Weyl group with the group $W$ of permutations $w\in S^{2N}$
with the property $w_i+w_{2N+1-i}=2N+1$ for all $i$. 

The group $W \subset S^{2N}$ is generated by the simple transposition
$s^C_N=(N,N+1)$ and by the products of two simple
transpositions $s_i^C = (i,i+1) (2N-i,2N-i+1)$, $i = 1, \dots , N-1$.

For a permutation $w$, denote 
$FL^\perp (V_P)  \cap  G_w^{{\mc F}^\infty}$ by $G_w^C$. 
We have $G_w^C \neq \emptyset$ if and
only if $w \in W$. Therefore we have 
a decomposition 
\bean\label{dec sp}
FL^\perp(V_P)\ =\ \sqcup_{w\in W} \ G_w^C .
\eean
The variety $FL^\perp(V_P)$ is isomorphic to the
 $SP_{2N}$ flag variety. Under this isomorphism decomposition
\Ref{dec sp} becomes the Bruhat cell decomposition.
The Bruhat cells $G_w^C$ are topological
cells, see \cite{H}.

Consider the generating morphism $\beta: FL(V_P) \to \bs P(\C[x])^{2N-1}$.
For a tuple $\bs y \in \bs P(\C[x])^{2N-1}$, denote
$\bs l(\bs y)=(\deg y_1,\dots, \deg y_{2N-1})$.
Two flags $\mc F_1$ and $\mc F_2$ belong to the same Bruhat cell in $FL(V_P)$ 
if and only if $\bs l(\beta(\mc F_1))=\bs l(\beta(\mc F_2))$.
Set as before $\bs l^w=\bs l(\beta(\mc F))$ for $\mc F \in G_w^{\mc F^\infty}$.
For $w\in W$ we have $l^w_i=l_{2N-i}^w$ for all $i$.

The Weyl groups of Langlands dual Kac-Moody algebras are naturally
identified. So the Weyl groups of $sp_{2N}$ and $so_{2N+1}$ are naturally
identified. Therefore we can consider an element $w\in W$ as an element of
the $so_{2N+1}$ Weyl group as well, in particular $w$ acts on  $so_{2N+1}$-weights.

\begin{lem} For $w\in W$, we have the equality of $so_{2N+1}$-weights:
\be
\sum_{i=1}^{N}\  l_i^w\, \al_i^B\ =\ \sum_{s=1}^n \ {\La_s^B}\ -\ w \cdot \tilde\La_\infty^B .
\ee
\end{lem}
\begin{proof} Follows from Lemma \ref{sln bruhat lem}. 
\end{proof}

Let $\bs y$ represent a critical point of the master function
$\Phi(\bs t; \bs z; \bs \La^B, w\cdot\tilde\La_\infty^B)$ and let $P$
be the $so_{2N+1}$-population originated at $\bs y$.

\begin{cor} The closure in $\C^0_{l_1^w}[x]\otimes\dots\C^0_{l_N^w}[x]$ of
  $N$-tuples of polynomials in the $so_{2N+1}$-population $P$
  representing a critical point of the master function $\Phi(\bs t; \bs z; \bs
  \La^B, w\cdot\tilde\La_\infty^B)$, is isomorphic to the Bruhat cell
  $G_w^C$ of the flag variety $FL^\perp (V_P)$.
\end{cor}

\subsection{The $sp_{2N}$-populations and Bruhat cells}
Let $P$ be an $sp_{2N}$-population with fundamental space $V_P$.
The space $V_P$ is a  selfdual space of dimension $2N+1$.
The canonical bilinear form of $V_P$ is symmetric.
Let $z_1, \dots , z_n, \infty$ be ramification points of $V_P$.
Let $\La_1^A, \dots , \La_n^A, \tilde\La_\infty^A$
be the corresponding integral dominant $sl_{2N+1}$-weights and 
$\La_1^C, \dots , \La_n^C, \tilde \La_\infty^C$ the corresponding integral dominant 
$sp_{2N}$-weights of $P$.

The generating morphism $\beta$ identifies the variety of isotropic flags
$FL^\perp (V_P)$ with the population $P$.
The variety $FL^\perp (V_P)$ is isomorphic to the flag variety of 
the Lie group $SO_{2N+1}$. 
The root system of the corresponding Lie algebra $so_{2N+1}$ is of type $B_N$.

Now we describe a Bruhat cell decomposition of the flag variety of $SO_{2N+1}$
as  a decomposition of $FL^\perp(V_P)$.

Let $\iota: \h^*_{so_{2N+1}}\to \h^*_{sl_{2N+1}}$ be the linear inclusion 
given by 
\be
\al_i^B \ \mapsto \ \al_i^A+\al_{2N-i}^A , 
\qquad \al_N^B\ \mapsto \ \al_N + \al_{N+1} ,
\ee
$i = 1, \dots , N-1$.  We have 
$(\al_i^B , \al_j^B)\ =\ (\iota(\al_i^B) ,  \iota(\al_j^B))$ and therefore 
$\iota$ is orthogonal. We have the following relation between simple
reflections:
\be
\iota (s_i^B \La)\ =\ s_i^A s^A_{2N-i} \iota(\La) ,
\qquad 
\iota (s_N^C\La)\ =\ s_N^A s_{N+1}^A s_N^A \iota(\La) ,
\ee
$i = 1, \dots , N-1$, for any weight $\La \in \h^*_{so_{2N+1}}$.
Therefore the map from the $so_{2N+1}$ Weyl group to the $sl_{2N+1}$
Weyl  group given by 
\be
s_i^B\ \mapsto \ s^A_i s^A_{2N-i} ,
\qquad 
s_N\ \mapsto \ s_N s_{N+1} s_N ,
\ee
$i = 1, \dots , N-1$ is a group homomorphism. The map identifies the
$so_{2N+1}$ Weyl group with the group $W$ of permutations $w \in S^{2N+1}$
with the property $w_i + w_{2N+2-i} = 2N + 2$ for all $i$.

The group $W \subset S^{2N+1}$ is generated by the simple
transposition  $s^B_N \, =\, (N,N+2)$ and by the products of two simple
transpositions $s_i^B \, =\, (i,i+1)(2N-i,2N-i+1)$, $i = 1, \dots , N-1$.

Let
\be
FL(V_P)\ =\ \sqcup_{w\in S^{2N+1}}\ G_w^{{\mc F}^\infty}
\ee
be the Bruhat cell decomposition with respect to the flag ${\mc F}^\infty$. 
For a permutation $w$, denote $FL(V_P)^\perp \cap G_w^{{\mc F}^\infty}$ by $G_w^B$. 
We have $G_w^B \neq \emptyset$ if and
only if $w \in W$. Therefore we have a decomposition
\bean\label{dec so}
FL^\perp(V_P)\ =\ \sqcup_{w\in W}\ G_w^B .
\eean
The variety $FL^\perp(V_P)$ is isomorphic to the $SO_{2N+1}$ flag variety.
Under this isomorphism decomposition \Ref{dec so}  becomes the Bruhat cell decompostion.

Consider the generating morphism $\beta : FL(V_P) \to \bs P(\C[x])^{2N}$.
For a tuple  $\bs y \in \bs P(\C[x])^{2N}$ denote
$\bs l(\bs y) = (\deg y_1, \dots , \deg y_{2N})$.
Two flags $\mc F_1$ and
$\mc F_2$ belong to the same Bruhat cell if and only if $\bs
l(\beta(\mc F_1)) = \bs l(\beta(\mc F_2))$.
Set as before $\bs l^w = \bs l(\beta(\mc F))$ for $\mc F\in G_w^{\mc F^\infty}$.
For $w \in W$ we have $l^w_i = l_{2N+1-i}^w$ for all $i$..

The Weyl groups of Langlands dual Kac-Moody algebras are naturally
identified. So the Weyl groups of $so_{2N+1}$ and $sp_{2N}$   are naturally
identified. Therefore we can consider an element $w \in W$ as an element of
the $sp_{2N}$ Weyl group as well, in particular $w$ acts on  $sp_{2N}$-weights.

\begin{lem} For $w \in W$, we have the equality of $sp_{2N}$-weights:
\be
\sum_{i=1}^{N} \ l_i^w \al_i^C\ =\ \sum_{s=1}^n {\La_s^C}\ -\ w\cdot \tilde\La_\infty^C .
\ee
\end{lem}
\begin{proof} Follows from Lemma \ref{sln bruhat lem}. 
\end{proof}

Let $\bs y$ represent a critical point of the master function
$\Phi(\bs t; \bs z; \bs \La^C, w\cdot\tilde\La_\infty^C)$ and let $P$
be the $sp_{2N}$-population originated at $\bs y$.

\begin{cor} The closure in $\C_{l_1^w}^0[x]\otimes \dots\otimes
  \C_{l_N^w}^0[x]$ of $N$-tuples of polynomials 
in the $sp_{2N}$-population $P$ representing critical points of the
master function
$\Phi(\bs t; \bs z; \bs  \La^C, w\cdot\tilde\La_\infty^C)$ is 
isomorphic to the Bruhat cell
  $G_w^B$ of the flag variety $FL^\perp(V_P)$.
\end{cor}

\section{Local Version}\label{sec loc}

The proof of  Theorems \ref{conv so} and \ref{conv sp} require a local variant of
the theory of populations where polynomials are replaced by germs of holomorphic functions.
Below we sketch this local variant.

\subsection{Generic and fertile tuples of  germs}\label{fertile germs} 
Let $z \in \C$. Denote $\C((x-z))$ the complex vector space of germs at $z$ of 
holomorphic functions in $x$. Denote $\bs P(\C((x-z)))$ the corresponding projective space.

Let $A$ be a  generalized Cartan $r \times r$-matrix  and
$\g(A)$ the corresponding Kac-Moody algebra. We use the same notations as in Section \ref{Kac_Moody sec}.

Fix an $r$-tuple $\bs T = (T_1, \dots , T_r) \in \Pzr$.
We say that $\La \in \h^*$  is {\it the weight of} $\bs T$ if for any $i$ the number
$\langle \La, \al_i^\vee\rangle$ is equal to the order of zero at $z$ of the germ $T_i$.
The weight $\La$ is integral dominant. 

For $\bs y = (y_1, \dots , y_r)\ \in \Pzr$, let $l_i$ be the order of zero at $z$ of $y_i$. The weight
$\La - \sum_{i=1}^r l_i \al_i$ is called {\it the weight of $\bs y$ with respect to $  \bs T,\ \La$}.

We say that an $r$-tuple  $\bs y  \in \Pzr$ is {\it generic} at $z$ 
if $y_i(z) \neq 0$ for all $i$. The weight of a generic tuple with respect to $ \bs T,\ \La$
is equal to $\La$.

We say that an $r$-tuple  $\bs y  \in \Pzr$ is {\it fertile}
with respect to $\bs T$ 
if for every $i$ there exists a germ $\tilde y_i$ satisfying equation
\Ref{wronskian-critical eqn}. 
If an $r$-tuple $\bs y$  is generic, then it is fertile. 

Let $\bs y$ be generic and $i = 1, \dots , r$. Lemma \ref{tilde y} shows 
that the $r$-tuples 
$\bs y^{(i)}$ $ = $ $(y_1, \dots , $ $ \tilde y_i,\dots, y_r) \ \in  \Pzr$,
where $\tilde y_i$ is given  by \Ref{tilde}, form a one-parameter family.
The parameter space of the family is identified with the projective line 
$\bs P^1$ with projective coordinates $(c_1 : c_2)$.
We have a map
\bean\label{map hol}
Y_{\bs y, i}\ : \ \bs P^1 \ \to \Pzr\ , 
\eean
which sends a point $c = (c_1 : c_2)$ to the corresponding $r$-tuple $\bs y^{(i)}$.
Almost all $r$-tuples $\bs y ^{(i)}$ are generic.
The exceptions form a finite set in $\bs P^1$.

\subsection{Reproduction procedure}\label{procedure for germs}
Let $\bs T$ be as in Section \ref{fertile germs}. Let $\bs y^0$ be generic.
Let $\bs i = (i_1, i_r, \dots , i_k), \ 1 \leq i_j \leq r,$ be a sequence of natural numbers.
We define a $k$-parameter family of fertile $r$-tuples 
\bea\label{general map germs}
Y_{\bs y^0, \bs i} \ :\ (\bs P^1)^k \ \to \Pzr
\eea
as in Section \ref{general procedure}. Namely, we proceed 
by induction on $k$,  starting at $\bs y^0$ and successively applying 
the simple reproduction procedure in directions
$i_1, \dots , i_k$.
The union 
\be 
P_{\bs y^0} \ = \ \cup_{\bs i} \ P_{\bs y^0, \bs i}\ \subset \Pzr\ ,
\ee
where the summation is over all sequences $\bs i$, 
is called {\it the population of germs at $z$ associated} with the Kac-Moody algebra $\g$, 
the tuple $\bs T$ of weight $\La$,   and {\it originated} at $y^0$.
We say that $\bs T$ is {\it the defining tuple} of the population $P$.

For a given $\bs i = (i_1, \dots , i_k)$, almost all $r$-tuples 
$Y_{\bs y^0, \bs i} ( \bs c)$ are generic. 
Exceptional values of $c \in (\bs P^1)^k$ are contained in a proper algebraic subset.

If two populations of germs at $z$ with the same $\bs T$ intersect, then they coincide.

As in Section \ref{Conjec} 
we conjecture that every population of germs at a point $z$, 
associated to a Kac-Moody algebra $\g$ and a tuple $\bs T$
is a (pro-)algebraic variety isomorphic to 
the full flag variety associated to the Kac-Moody algebra
$\g^t$ which is Langlands dual to $\g$. Moreover, the parts of the family corresponding to 
$r$-tuples of germs with fixed degrees of zero at $z$ are isomorphic to 
Bruhat cells of the flag variety.

\subsection{Orders of zero of germs in a population and the Weyl group}
 Let a tuple $\bs y$ be a member of a population $P$
associated to a tuple $\bs T$ of weight $\La$.
Let a tuple $ {\bs y}^{(i)} = (y_1, \dots , \tilde y_i, \dots , y_r)$ be 
an immediate descendant of $\bs y$.  
Let  $ \La^{(i)}$ be the weight of $ {\bs y}^{(i)}$ with respect to $\bs T, \La$.
If the order of zero at $z$ of $\tilde y_i$ is not equal to the order of zero of  $y_i$, then  
$\Lambda^{(i)}_\infty\ = \ s_i\cdot{\La}_\infty$, cf. Lemma 
\ref{simple inf lem}.

Let $P$ be the population of germs associated with a tuple $\bs T$ of weight
$\La$ and originated at a tuple $\bs y^0$. Then
for any  $\bs y \in P$,  there is an element $w \in \mathcal W$, 
such that the weight of $\bs y$ is  $w \cdot \Lambda$, and for
any  $w \in  \mathcal W$,  there is a tuple $\bs y \in P$ with weight
$w \cdot \Lambda$.

\subsection{Relation between populations of polynomials and germs}
Let $P$ be a population of critical points in the sense of Section
\ref{general procedure}. We assume that $P$ is associated with integral dominant weights
$\La_1, \dots \La_n$, complex numbers $z_1, \dots , z_n$, and originated at a critical
point $\bs y^0$. Let $\bs T = (T_1, \dots , T_r)$ be the tuple of polynomials associated with $P$ and 
defined by formulas \Ref{T}.

Let $z \in \C$. Then $P$ induces a population of germs at $z$, denoted $P_z$. 
The population of germs $P_z$ consists of germs at $z$ of tuples $\bs y \in P$. 
The tuple $\bs T$, considered as an element of $\Pzr$, is the defining tuple of $P_z$.
The weight of $\bs T$ at $z$ is $\La_i$ if $z = z_i$ and is zero otherwise.

In order to check that $P_z$ is indeed a population of germs  one needs to check that
$P_z$ contains a tuple of germs $\bs y$ such that $y_i(z) \neq 0$ for all $i$.
To do this we consider the tuple of polynomials $\bs y^0$. The tuple is generic
in the sense of Section \ref{PLCP}. Thus if $z$ is not a root of one of $y_i$,
 then the germ of $\bs y^0$ is the required tuple. If $z$ is a root of one of $y_i$,
then it is a simple root. It is easy to see that any immediate descendant of $\bs y^0$
in the $i$-th direction gives a required tuple.

\subsection{$sl_{N+1}$-Populations of germs}\label{germs sl}
Let $P$ be an $sl_{N+1}$-population of germs at $z \in \C$ with defining
tuple $\bs T$ of weight $\La$. Let $\bs y \in P$. Introduce a germ of a 
linear differential operator $D(\bs y)$ by formula \Ref{Diff}. The germ $D(\bs y)$
does not depend on the choice of $\bs y$. It is called associated to $P$ and denoted
$D_P$. The exponents of $D_V$ are the numbers
$0, \ (\La + \rho, \al_1),\  \dots , \ (\La + \rho, \al_1 + \dots + \al_N)$.
The kernel $V$ of $D_V$ is called the fundamental space of $P$ and denoted $V_P$.

Let $V$ be a complex vector space of germs at $z \in \C$ of holomorphic functions
in $x$, $\dim V = N+1$. Assume that $V$ does not have a base point, i.e. 
there is $f \in V$ such that $f(z) \neq 0$. Let $D_V = \partial^{N+1} + ...$
be the germ of a linear differential operator with kernel $V$. Let 
$0,\ m_1 + 1, \ \dots , \ m_1 + \dots + m_N + N$ be the exponents at $z$ of $D_V$.
Here $m_i$ are non-negative integers. Introduce an integral dominant $sl_{N+1}$-weight
$\La$ by  conditions $(\La, \al_i) = m_i$.

Let $u_1, \dots , u_{N+1}$ be a basis in $V$. The germ 
$W_V = W(u_1, \dots , u_{N+1})$ is uniquely  determined up to multiplication
by a non-zero number. Fix a tuple $\bs T = (T_1, \dots , T_N) \in \PzN$ so that
\begin{enumerate}
\item[$\bullet$]
The product $T_1^{N} T_2^{N-1} \dots T_N$ is equal to $W_V$ up to 
multiplication by a nonzero number;
\item[$\bullet$]
For any $i$, the order of zero at $z$ of $T_i$ is $m_i$.
\end{enumerate}
The $sl_{N+1}$-weight of $\bs T$ is $\La$.
Such a tuple $\bs T$ is called {\it a framing} of $V$.

Let $FL(V)$ be the variety of full flags in $V$. Define the generating morphism
\be
\beta \ :\ FL(V) \ \to \ \PzN , 
\qquad
\mc F \ \mapsto \ \bs y^{\mc F} ,
\ee
by the same formulas as in Section \ref{ha-ha}. We have the following analog of 
Theorem \ref{pol-crit}.

\begin{theorem}\label{A germs}

${}$

\begin{enumerate}
\item[1.] The generating morphirm  defines an isomorphism of \ $FL(V)$
and its image \newline $\beta (FL(V)) \subset \PzN$.
\item[2.]
The image $\beta(FL(V))$ is an $sl_{N+1}$-population of germs associated with
$\bs T$ and $\La$.
\item[3.]
The initial vector space $V$ is the fundamental space of the population $\beta (FL(V))$.
\item[4.]
Every $sl_{N+1}$-population of germs is given by this construction
starting from a suitable space $V$.
\end{enumerate}
\end{theorem}

\subsection{Selfdual vector spaces of germs}
Let $V$ be an $N+1$-dimensional vector space of germs at $z \in \C$ of holomorphic functions.
For $u_1, \dots , u_i \in V$ define
$W^\dagger (u_1, \dots , u_i)$, the divided Wronskian with
with respect to $V, \bs T$, by the same formulas as in Section \ref{DU}.
Define the vector space $V^\dagger$ as the span of divided Wronskians
$W^\dagger(u_1, \dots , u_N)$ with $u_1, \dots u_N \in V$.

An $N$-tuple $\bs y \in \PzN$ is called {\it symmetric} if for any $i$
the germ $y_i$ is equal to the germ $y_{N+1-i}$ up to multiplication 
by a non-zero number.

The space $V$ with a symmetric framing $\bs T$ is called {\it selfdual} if 
$V^\dagger = V$.

A selfdual space has a nondegenerate bilinear form defined
as in Section \ref{SELF}. The form is skew-symmetric if the dimension of
$V$ is even and is symmetric if the dimension is odd.

We denote $FL^\perp(V) \subset FL(V)$ the subvariety of isotropic flags.
Let $\beta$ be the generating morphism of $V$. A flag $\mc F$ is isotropic
if and only if the tuple $\beta (\mc F) = \bs y^{\mc F}$ is symmetric.
The group of linear transformations of $V$ preserving the
bilinear form acts on $\beta (FL^\perp(V))$ be the same formulas as in Section \ref{gen iso}.

\subsection{Populations of germs of type $B_N$} We use notations of Section \ref{BN}.
Let $P$ be an $so_{2N+1}$-population of germs associated with a tuple 
$\bs T = (T_1, \dots , T_N)$ of weight $\La$
and originated at a generic $N$-tuple $\bs y^0$. 

For any $N$-tuple $\bs y$ denote the symmetric $2N-1$-tuple
$(y_1, \dots , y_{N-1}, y_N, y_{N-1}, \dots , y_1)$  by $\bs y_A$.
The $2N-1$-tuple $\bs T_A$ has the $sl_{2N}$-weight $\La^A$. 

Let $P^A$ be the $sl_{2N}$-population of germs associated with
the tuple $\bs T^A$ of $sl_{2N}$-weight $\La^A$ and originated at 
$\bs y^0_A$. There is an injective map $P \to P^A, \ \bs y \mapsto \bs y_A$.
The fundamental space of $P_A$ is selfdual with respect to the symmetric
framing $\bs T_A$. 
The fundamental space of of $P_A$ is called the fundamental space of $P$.
The generating morphism $\beta  : FL(V)  \to  \PzN $ identifies the subvariety $FL^\perp(V)$
of flags in $V$, isotropic with respect to its skew-symmetric form, and the 
$so_{2N+1}$-population $P$.

\begin{theorem}\label{loc BN}
Let $V$ be a vector space of germs of dimension $2N$ with no base point. Assume that $V$ is selfdual
with respect to a symmetric framing $\bs T_A$. Let $FL^\perp(V)$ be the variety of isotropic 
full flags of $V$, and $\beta(FL^\perp(V))$ its image under the generating morphism. Then
$\beta(FL^\perp(V))$ is an $so_{2N+1}$-population of germs with fundamental space $V$.
\end{theorem}
To prove this theorem it is enough to show only that there is an isotropic flag $\mc  F$ such that its image
$y^{\mc F}$ is generic, i.e. $y^{\mc F}_i(z) \neq 0$ for all $i$. This is proved similarly to the proof of
Theorem \ref{special basis prop}.

\subsection{Populations of germs of type $C_N$} We use notations of Section \ref{CN}.
Let $P$ be an $sp_{2N}$-population of germs associated with a tuple 
$\bs T = (T_1, \dots , T_N)$ of weight $\La$
and originated at a generic $N$-tuple $\bs y^0$. 

Denote $\bs {\tilde T}$ the symmetric $2N$-tuple
$(T_1, \dots , T_{N}, T_N, \dots , T_1)$.
The tuple $\bs {\tilde T}$ has the $sl_{2N+1}$-weight $\La^A$. 

For any $N$-tuple $\bs y$ denote $\bs y_A$ the $2N$-tuple
$(y_1, \dots , y_{N-1}, y_N^2, y_N^2, y_{N-1}, \dots , y_1)$.

Let $P^A$ be the $sl_{2N+1}$-population of germs associated with
the symmetric tuple $\bs {\tilde T}$ of $sl_{2N+1}$-weight $\La^A$ and originated at 
$\bs y^0_A$. There is an injective map $P \to P^A, \ \bs y \mapsto \bs y_A$.
The fundamental space of $P_A$ is selfdual with respect to the symmetric
framing $\bs {\tilde T}$.
The fundamental space of $P_A$ is called the fundamental space of $P$.
The generating morphism $\beta  : FL(V)  \to  \PzN $ identifies the subvariety $FL^\perp(V)$
of flags in $V$, isotropic with respect to its symmetric form, and the 
$so_{2N+1}$-population $P$.

\begin{theorem}\label{loc CN}
Let $V$ be a vector space of germs of dimension $2N+1$ with no base point. Assume that $V$ is selfdual
with respect to a symmetric framing $\bs {\tilde T}$. Let $FL^\perp(V)$ be the variety of isotropic 
full flags of $V$, and $\beta(FL^\perp(V))$ its image under the generating morphism. Then
$\beta(FL^\perp(V))$ is an $sp_{2N}$-population of germs with fundamental space $V$.
\hfill $\square$
\end{theorem}

\subsection{Proof of Theorems \ref{conv so} and \ref{conv sp}} We prove Theorem \ref{conv so}.
The proof of Theorem \ref{conv sp} is similar.

Let $V$ be a vector space of polynomials of dimension $2N$ selfdual in the sence of Section
\ref{SELF}. Let $z \in \C$. Denote $\pi_z$ the map which sends a polynomial to its germ at $z$.
Then $\pi_z(V)$ is a $2N$-dimensional vector of germs.

The symmetric polynomial framing  
$\bs T = (T_1, \dots , T_{2N-1})$, defined for $V$ by formulas \Ref{T polyn}, induces a symmetric framing
$\pi_z(\bs T)$ of $\pi_z(V)$. The vector space $\pi_z(V)$ is selfdual with respect to the framing
$\pi_z(\bs T)$. The bilinear forms on $V$ and $\pi_z(V)$ are functorial, for any $f, g \in V$ we have
$(f, g) = (\pi_z(f), \pi_z(g))$. The map $\pi_z$ identifies 
the varieties of isotropic flags $FL^\perp(V)$ and $FL^\perp(\pi_z(V))$.

The map $\pi_z$ commutes with the generating morphism $\beta$.

According to Theorem \ref{loc BN} the image $\beta(FL(\pi_z(V)))$ is an $so_{2N+1}$-population
of germs. In particular, $\beta(FL^\perp (\pi_z(V)))$ contains $2N-1$-tuples generic
in the sense of Section
\ref {fertile germs}. Our goal is to show that $\beta(FL^\perp (V))$ contains a tuple of polynomials
generic in the sense of Section \ref{PLCP}. 
We start with an arbitrary flag $\mc F^0 \in FL^\perp(V)$.
We  will deform it slightly so that the corresponding tuple will
become generic. Let $\beta (\mc F^0) = \bs y^0$ and $\bs y^0 = (y_1^0, \dots , y_{N-1}^0,
$ $y_N^0,$ $  y_{N-1}^0,$ $ \dots , y^0_1)$. We may assume that 
the roots of  polynomials $y^0_i$ do not contain the ramification points $z_1, \dots , z_n$ of $V$.
Indeed by Theorem \ref{loc BN}  we know that
for every $z$ almost all flags in $FL^\perp(\pi_z(V))$ give generic tuples. 
The set of ramification points is finite. So almost all flags
in $FL^\perp(V)$ give generic tuples of germs at the ramification points.
Thus we may assume that
if $z \in \{z_1, \dots , z_n\}$, then the tuple of germs $\pi_z(\bs y^0)$ is generic.

Let $Z(\bs y^0) \subset \C$ be the finite set of all $z$ such that $\pi_z(\bs y^0)$ is not generic.
If $\bs y^0$ is slightly deformed in $\beta(FL^\perp (V))$, then $Z(\bs y^0)$ is slightly deformed in $\C$.
In particular the slightly deformed
$Z(\bs y^0)$ still will not intersect  the set of ramification points of $V$.
Our goal is to show  that there is  a small deformation $\bs y^1$ of $\bs y^0$ in $\beta(FL^\perp(V))$
such that for any $z \in Z(\bs y^1)$
exactly one of polynomials $y^1_1, \dots , y^1_N$ has a root at $z$ and the multiplicity of that
 root is equal to one.

Assume that $\bs y^0$ is slightly deformed in the direction of one of the one-parameter groups indicated in
Lemmas \ref{skew 1} and \ref{skew 2}. Let $\bs y^1$ be the deformed tuple. 
If $z$ belongs to $Z(\bs y^1)$ and $z$ does not belong to $Z(\bs y^0)$, then 
exactly one of the polynomials  $y_1^1, \dots , y^1_N$ has a root at $z$ and the multiplicity of
that root is one.

Thus our goal is to deform   slightly the tuple $\bs y^0$ 
in the directions of the  one-parameter subgroups of Lemmas 
\ref{skew 1} and \ref{skew 2}, so that for the deformed tuple $\bs y^1$ the sets $Z(\bs y^0)$ and
$Z(\bs y^1)$ do not intersect. 

The possibility to do so is based on the following important observation. For any given $z$, a 
small generic deformation of $\bs y^0$ in the direction of one of the one-parameter subgroups of Lemmas 
\ref{skew 1} and \ref{skew 2} does not increase the order of zeros of germs of $\pi_z(\bs y^0)$.

Let $z_0$ be one of points in $ Z(\bs y^0)$. By Theorem \ref{loc BN}  we know that
 $\pi_{z_0}(\bs y^0)$ can be made generic using small
deformations in the directions of the distinguished one-parameter subgroups. According to the above observation
we will not make the germs $\pi_z(\bs y^0),\ z \,\in \,Z(\bs y^0) - \{z_0\}$, worse while simplifying the germ at $z_0$.

Those remarks prove Theorem \ref{conv so}.

\section{Appendix: the Wronskian identities}
In this appendix we collect identities involving Wronskians.
All functions in this section are functions of $x$ with sufficiently
many derivatives.

Recall that the Wronskian of functions $g_1,\dots,g_s$ is defined by
\be
W(g_1,\dots,g_s)\ =\ \det(g_i^{(j-1)})_{i,j=1}^s\ .
\ee
We follow the convention that for $s=0$ the corresponding Wronskian equals $1$. 
In this section we write $W_s(g_1,\dots,g_s)$ instead of $W(g_1,\dots,g_s)$ to
stress the order of the Wronskian.

We start with
\begin{lem}\label{1 wronskian}
We have $W_{s+1}(1,g_1,\dots,g_s)\ =
\ W_s(g_1',\dots,g_s')$.
\end{lem}
\begin{proof} Obvious. \end{proof} 

The next lemma is

\begin{lem}\label{f wronskian}
We have $ W_s(fg_1,\dots,fg_s)\ =\ f^s\ W_s(g_1,\dots,g_s)$.
\end{lem}
\begin{proof}
We proceed by induction on $s$. The case $s=0$ is obvious. Suppose the
lemma is proved for $s=s_0-1$. Compare the following 
differential equations on $g_1$
\be 
W_{s_0}(fg_1,\dots,fg_{s_0})\ =\ 0\  \qquad 
\on{and} \qquad
W_{s_0}(g_1,\dots,g_{s_0})\ =\ 0\ .
\ee
Both equations have solutions $g_2,\dots,g_{s_0}$
and therefore coincide up to multiplication by a function. Comparing
the coefficients of $g_1^{(s_0-1)}$ we find that this function is 
$f^{s_0}$ by the induction hypothesis.
\end{proof}

In what follows we 
use Lemmas \ref{1 wronskian} and \ref{f
  wronskian} to prove new identities.

\begin{lem}\label{fg id} We have 
\be
W_{s+1}(f^s,f^{s-1}g,f^{s-2}g^2,\dots,g^s)\ = 
 \ (\prod_{i=1}^s i!) \  W_2(f,g)^{s(s+1)/2}\ .
\ee
\end{lem}
\begin{proof}
Let $h=g/f$. We have
\be
W_{s+1}(f^s,f^{s-1}g,f^{s-2}g^2,\dots,g^s)\ =\ f^{s(s+1)}\ W_{s+1}(1,h,h^2,\dots,h^s)\ .
\ee
This Wronskian is equal to the Wronskian of order $s$
\be
f^{s(s+1)}\ W_s(h',2hh',\dots,sh^{s-1}h')\ 
=\ s!\ (h')^s\ f^{s(s+1)}\ W_s(1,h,h^2,\dots,h^{s-1})\ .
\ee
Continuing this reduction, we obtain
\be
W_{s+1}(f^s,f^{s-1}g,f^{s-2}g^2,\dots,g^s)\ =\ (\prod_{i=1}^s i!)\ f^{s(s+1)}\ (h')^{s(s+1)/2}\ .
\ee
On the other hand
\be
W_2(f,g)^{s(s+1)/2}\ =\ f^{s(s+1)}\ W(1,h)^{s(s+1)/2}\ =\ f^{s(s+1)}\ (h')^{s(s+1)/2} \ ,
\ee
and the lemma follows.
\end{proof}

For the next two identities, we 
fix integers $0\leq k\leq s$ and functions $g_1,\dots,g_{s+1}$. 
Let
\be
V_{s-k+1}(i)\ =\ W_{s-k+1}(g_1,\dots,g_{s-k}, g_i).
\ee
\begin{lem}\label{wr id 2} We have
\be
W_{k+1}(V_{s-k+1}(s-k+1),\dots, V_{s-k+1}(s+1))=
(W_{s-k}(g_1,\dots,g_{s-k}))^k W_{s+1}(g_1,\dots,g_{s+1}) .
\ee
\end{lem}
\begin{proof}
This lemma is proved by induction on $s$. The case $k=s$ is
trivial. Suppose that the lemma is proved
for $s=k,\dots,s_0-1$. 

Divide both sides of the identity for $s=s_0$
by $g_1^{(s_0-k+1)(k+1)}$ and use Lemma \ref{f wronskian} to carry $g_1$
inside  all Wronskians. Then one of the functions in each Wronskian
is 1 and we can reduce the order by Lemma \ref{1 wronskian}.
Then the identity for $s=s_0$ follows from the  induction 
hypothesis applied to $f_i = (g_{i+1}/g_1)'$, $i = 1, \dots, s_0$.
\end{proof}

Let
\be
W_s(i)\ =\ W(g_1,\dots,\widehat{g_i},\dots,g_{s+1})
\ee
be the Wronskian of all functions except  $g_i$.
\begin{lem}\label{wr id 1}
We have
\be
W_{k+1}(W_s(s+1),W_s(s),\dots,
W_s(s-k+1))=W_{s-k}(g_1,\dots,g_{s-k})(W_{s+1}(g_1,\dots,g_{s+1}))^{k}.
\ee
\end{lem}
\begin{proof}
First we prove the case $s=k$ by induction on $k$. The case $k=0$ is
trivial. Suppose the case $k<k_0$ is proved. Divide both sides of our
identity in the case $s=k=k_0$ by $g_1^{k_0(k_0+1)}$. By Lemmas
\ref{1 wronskian} and \ref{f wronskian}
we are reduced to the identity
\be
W_{k_0+1}(W^{h'}_{k_0-1}(k_0),\dots,W^{h'}_{k_0-1}(1),W_{k_0}(h_1,\dots,h_{k_0}))
\ =\ (W(h_1',\dots,h_{k_0'}))^{k_0}\ ,
\ee
where $h_i = g_{i+1}/g_1$ and $W^{h'}_{k_0-1}(i) = W_{k_0-1}(h_1',\dots,\widehat{h_i'},\dots,h_{k_0}')$.

The left hand side of the last identity is a determinant of size $k_0+1$. Add to the last
row the row number $i$ with coefficient $(-1)^{k_0-i+1}h_{k_0-i+1}$,
$i=1,\dots,k_0$. Then the last row becomes
\be
(0,\dots,0,W_{k_0}(h_1',\dots,h_{k_0}'))
\ee
and the lemma for $k=k_0$ follows from the induction hypothesis
applied to functions $f_1=h_1',\dots,f_{k_0}=h_{k_0}'$.

Now we continue by induction on $s$. Suppose that the lemma is proved
for $s=k,\dots,s_0-1$. Divide both sides of the identity for $s=s_0$
by $g_1^{s_0(k+1)}$. Then the identity for $s=s_0$ follows from the  induction
  hypothesis applied to $f_i=(g_{i+1}/g_1)'$, $i=0,\dots,s_0-1$.
\end{proof}


\begin{thebibliography}
\normalsize


\bibitem[EH]{EH} D. Eisenbud, J. Harris, {\it Limit Linear Series: Basic Theory},
Inventiones mathematicae, 85, 337-371.



\bibitem[Fa]{Fa} L. Faddeev, Lectures on the Quantum Inverse Scattering Method.
In: Integrable Systems, ed. by X.-G. Song, Nankai Lectures Math. Phys., 
World Scientific, Teaneck, N.J., 1990, 23-70.

\bibitem[FFR]{FFR} B. Feigin, E. Frenkel, and N. Reshetikhin, 
{\it Gaudin model, Bethe ansatz, and critical level},
CMP, 166 (1994), 27-62.

\bibitem[Fr]{Fr} E. Frenkel,  {\it Affine Algebras, Langlands Duality and Bethe Ansatz},
Proceedings of the International Congress of Math. Physics, Paris 1995,
International Press, 1995, 606-642,  q-alg/9506003.

\bibitem[Fu]{Fu} W. Fulton, {\it Intersection Theory}, Springer-Verlag, 1984.


\bibitem[G]{G} M. Gaudin, Diagonalization d'une class hamiltoniens de spin. 
{\it Journ. de Physique} {\bf 37}, no. 10 (1976), 1087 - 1098.


\bibitem[GH]{GH} Ph. Griffiths, J. Harris, {\it Principles of
    algebraic geometry}, A Whiley-Interscience Publication, 1994.


\bibitem[H]{H} S. Helgason, {\it Differential geometry, Lie groups and
    symmetric spaces}, Academic Press, New York, 1978.

\bibitem[K]{K} V. Kac, {\it Infinite-dimensional Lie algebras},
  Cambridge University Press, 1990.

\bibitem[M]{M} Y. Macdonald, {\it Symmetric Functions and Hall Polynomials},
Oxford Univ. Press, 1979


\bibitem[MV]{MV} E. Mukhin and A. Varchenko,
Remarks on Critical Points of Phase Functions and Norms of Bethe Vectors.
In: Arrangements - Tokyo 1998. 
{\it Advanced Studies in Pure Mathematics} {\bf 27} (2000), 239 - 246.


\bibitem[RV]{RV} N. Reshetikhin, A. Varchenko, {\it Quasiclassical
asymptotics of solutions to the KZ equations}, Geometry, topology
    $\&$ physics, Conf. Proc. Lecture Notes Geom. Topology, VI,
    Internat. Press, Cambridge, MA (1995), 293-322.



\bibitem[Sk1]{Sk1} E. Sklyanin, Separation of variables in the Gaudin model.
{\it J. Soviet Math.} {\bf 47} ( 1989 ), 2473 - 2488.


\bibitem[Sk2]{Sk2} E. Sklyanin, 
The functional Bethe ansatz. In: Integrable and Superintegrable Systems,
ed. by B. Kupershmidt, World Scientific, Singapore, 1990, 8 - 33.

\bibitem[S]{S} G. Szego, {\it Orthogonal polynomials}, AMS, 3d edition, 1967.

\bibitem[ScV]{ScV} I.Scherbak and A. Varchenko, {\it Critical point of
    functions, $sl_2$ representations and Fuchsian differential
    equations with only univalued solutions}, math. QA/0112269, (2001)
    1-25.

\bibitem[SV]{SV} V. Schechtman, A. Varchenko, {\it Arrangements of
    hyperplanes and Lie algebra homology}, Invent. Math., {\bf 106}
    (1991), 139-194.


\end{thebibliography}
\end{document}